\documentclass[10pt]{article}
\usepackage{amsmath,latexsym,amsfonts,amssymb,isolatin1,color}
\usepackage[all]{xy}
\newtheorem{theorem}{Theorem}[section]
\newtheorem{lemma}[theorem]{Lemma}
\newtheorem{proposition}[theorem]{Proposition}
\newtheorem{coro}[theorem]{Corollary}
\newtheorem{conjecture}[theorem]{Conjecture}
\begin{document}
\newcommand{\cw}{\color{white}}
\newcommand{\pn}{\mathbb P (E)} \newcommand{\pnd}{\mathbb P
  (E^*)} \newcommand{\pnu}{\mathbb P (\mathbb R^{n})}
\newcommand{\pnud}{\mathbb P ((\mathbb R^{n})^*)}
\newcommand{\sln}{PSL(n,\mathbb R)}
\newcommand{\sld}{PSL(2,\mathbb R)}
\newcommand{\qf}{\mathcal{QF}}
\newcommand{\vm}{\vert\mu\vert} \newcommand{\proof}{{\sc
    Proof: }} \newcommand{\qed}{{\sc Q.e.d.}}
\newcommand{\grf}{\pi_1 (S)}
\newcommand{\seq}[1]{ \{#1\}_{m\in\mathbb N}}
\newcommand{\mapping}[4] { \left\{
\begin{array}{rcl}
#1 &\rightarrow& #2\\
#3 &\mapsto& #4 .
\end{array}
\right.  } \newcommand{\auteur}{\vskip 2truecm
\centerline{François Labourie} \centerline{Topologie et
  Dynamique} \centerline{Université Paris-Sud}
\centerline{F-91405 Orsay (Cedex)} } \title{Anosov Flows,
  Surface Groups and Curves in Projective Space\\
} \author{François LABOURIE \thanks{L'auteur remercie
    l'Institut Universitaire de France.}}  \maketitle
\section{Introduction}
{In his beautiful paper \cite{H}, N. Hitchin studied the
connected components of the space
$$
Rep(\grf,PSL(n,\mathbb R))=Hom^{red}(\grf,PSL(n,\mathbb
R))/PSL(n,\mathbb R),
$$}
of reducible representations of the fundamental group of
a compact surface $S$ into $PSL(n,\mathbb R)$. By reducible,
we mean representations that act as sum of irreducible
representations on the Lie algebra by the adjoint
representation. Using Higgs bundles techniques, he proved
two remarkable results. The first one deals with the number
of components of this space.
\begin{theorem}\label{Hitch1}{\sc [Hitchin]}
  If $n>2$, the space $Rep(\grf,PSL(n,\mathbb R))$ has three
  connected components when $n$ is odd, and six when $n$ is
  even
\end{theorem}

Notice that W. Goldman gave a complete description of these
connected components in the case of finite covers of
$PSL(2,\mathbb R)$ in \cite{G1}. In the case of
$PSL(2,\mathbb R)$, two homeomorphic components, called
Teichmüller spaces, play a central role. These two
components are well known to be homeomorphic to a ball of
dimension $6g-6$.

N. Hitchin has generalised this situation to $PSL(n,\mathbb
R)$.  Indeed, one of these components when $n$ is odd, and
two when $n$ is even, have a very simple topology. Let's
define a {\em $n$-Fuchsian} representation to be a
representation $\rho$ which can be written as
$\rho=\iota\circ\rho_{0}$, where $\rho_{0}$ is a cocompact
representation with values in $PSL(2,\mathbb R)$ and $\iota$
is the irreducible representation of $PSL(2,\mathbb R)$ in
$PSL(n,\mathbb R)$, We denote by $Rep_{H}(\grf,PSL(n,\mathbb
R))$ a connected component that contains Fuchsian
representations, and call it a {\em Hitchin's component}.
Actually there is one Hitchin's component when $n$ is odd
and two isomorphic when $n$ is even. N. Hitchin proved in
\cite{H}

\begin{theorem}\label{Hitch2}{\sc [Hitchin]}
  Each component $Rep_{H}(\grf,PSL(n,\mathbb R))$ is
  homeomorphic to a ball of dimension $\chi(S)(1-n^{2})$.
\end{theorem}
This last result actually extends to the case of adjoint
groups of real split forms. Although Hitchin's proof gives
an explicit parametrisation of this component, the
construction by itself sheds no light on the geometry
underlying these representations. With Higgs bundle
techniques, one can prove the representations in Hitchin's
component are irreducible (Lemma \ref{corohitchin}), but it
seems hard to  detect with these techniques whether these representations are faithful,
discrete, as well as to  understand the action of the group of outer
automorphisms of $\grf$ acts properly on this specific
component.

Nevertheless, the geometric significance of this component
is well known in dimension 2 and 3. For $n=2$, it is
Teichmüller component, corresponding to holonomies of
hyperbolic structures on $S$. For $n=3$, S. Choi and W.
Goldman proved in \cite{CG}
\begin{theorem}\label{CG2}{\sc [Choi-Goldman]} For $n=3$,
  Hitchin's component con\-sists of ho\-lo\-no\-mies of
  convex real projective structures on $S$. That is, for
  every representation $\rho$ in $Rep_{H}(\grf,PSL(3,\mathbb
  R))$, there exists an open set $\Omega$ in $\mathbb P
  (\mathbb R^{3})$ such that $\Omega/\rho(\grf)$ is
  homeomorphic to $S$.
\end{theorem}
As a consequence of this result, a representation $\rho$ in
Hitchin's component when $n=3$ preserves a $C^{1}$-convex
curve in $\mathbb P(\mathbb R^{3})$, namely the boundary of the open set $\Omega$ obtained by the previous theorem.

Our first result generalises this last situation. Let's introduce a definition.  A curve $\xi$ with values in
$\mathbb P(\mathbb R^{n})$ is said to be {\em hyperconvex}
if for any distinct points $(x_{1},\ldots,x_{n})$ the
following sum is direct
$$
\xi(x_{1})+\ldots+\xi(x_{n}).
$$
Furthermore, we say a hyperconvex curve is a {\em Frenet
  curve}, if there exists a family of maps
$(\xi^{1},\xi^{2},\ldots,\xi^{n-1})$ with
$\xi^{p}\subset\xi^{p+1}$, called the {\em osculating flag of $\xi$},
such that \begin{itemize} \item $\xi=\xi^{1}$ and $\xi^{p}$
  takes values in the Grassmannian of $p$-planes, \item if
  $(n_1,\ldots,n_l)$ are positive integers such that $
  \sum_{i=1}^{i=l}n_i\leq n, $ if $(x_{1},\ldots,x_{l})$ are
  distinct points, then the following sum is direct
\begin{eqnarray}
\xi^{n_1}(x_1)+\ldots+\xi^{n_{l}}(x_{l});\label{fre2}
\end{eqnarray}
\item finally, for every $x$, let $p=n_{1}+\ldots+n_{l}$,
  then
\begin{eqnarray}
\lim_{\xymatrix@1{(y_1,\ldots,y_l)\ \  \ \ar[r]_-{ y_i \ \hbox{\tiny all distinct}}& \ \ x}}
(\bigoplus_{i=1}^{i=l}\xi^{n_i}(y_i))=\xi^{p}(x).\label{fre3}
\end{eqnarray}
\end{itemize}
One notices that for a Frenet hyperconvex curve, $\xi^{1}$
completely determines $\xi^{p}$. Also if $\xi^{1}$ is
$C^{\infty}$, then $\xi^{p}(x)$ is  generated by
the derivatives at $x$ of $\xi^{1}$ up to order $p-1$.
However, in general, a Frenet hyperconvex curve has no
reason to be $C^{\infty}$ also its image is obviously a
$C^{1}$-submanifold.  Our main result is the following
\begin{theorem}\label{mainA} For every representation $\rho$
  in Hitchin's component, there exists a $\rho$-equivariant
  hyperconvex Frenet curve from $\partial_{\infty}\grf$ to
  $\mathbb P(\mathbb R^{n})$.
\end{theorem} 
Notice that the Veronese embedding from $\mathbb P(\mathbb
R^{2})$ to $\mathbb P(\mathbb R^{n})$ is a $SL(2,\mathbb R)$-equivariant hyperconvex Frenet curve; it corresponds to Fuchsian
representations. Our main theorem therefore says that the Veronese embedding persists under eventually large deformation of the representation.

Recall that we say a representation $\rho$ of $\Gamma$ with
values in a semi-simple Lie group $G$ is {\em purely
  loxodromic}, if for every $\gamma$ in $\Gamma$ different
from the identity, $\rho(\gamma)$ is conjugate to an element
in the interior of the Weyl chamber. For $G=PSL(n,\mathbb
R)$, this just means that the eigenvalues of $\rho(\gamma)$
are real and  with multiplicity 1.  For $G=PSL(2,C)$, we recover the classical notion.  As a consequence of the
techniques involved in the proof, we also obtain
\begin{theorem}\label{mainB} Every representation in Hitchin's component is discrete, faithful and purely loxodromic.  \end{theorem}

This theorem is a generalisation of the classical result for Teichmüller Space. It bears also some relations with a beautiful recent result of M. Burger, A. Iozzi and  A. Wienhard \cite{BIW}, announced while  the second draft of this paper was completed. They proved in particular that surface group representations with maximal Toledo invariant also have discrete images. Although the methods and the target groups are different (so that the results have non empty intersection, but none contains the other), it appears after a discussion together that dynamical ideas quite similar to those appearing in this paper can be applied to their situation, improving some geometrical aspects. It is also quite surprising that two classes of groups, namely isometry groups of Hermitian symmetric spaces on one hand, and $SL(n,\mathbb R)$ (and quite plausibly all real split forms) on the other hand, have some common features, not shared for instance with $PSL(2,\mathbb C)$.

We shall also prove in a following  paper that the mapping class group acts properly on the Hitchin component \cite{FL4}

We shall also state converse or refinements of these results in Section \ref{state}. We now describe more precisely the structure of this paper, then proceed to discussion and conjectures.

I wish to thank Bill Goldman, Nigel Hitchin, Rick Kenyon and
Frédéric Paulin for their interest in the subject and
many discussions as well as David Fried for an enlightening
comment on structural stability, and Olivier Biquard for
help on Higgs bundles. Finally, Alessandra Iozzi and Marc Burger are warmly thanked for pointing out an error in the original version of this paper. 

\subsection{Description of the article}
We describe briefly the content of the main sections of this
article.
\begin{itemize}
\item \ref{1}: {\bf Geometric Anosov flows.} We introduce in
  this section a notion of ``geometric structure'', called
  {\em Anosov structure}, related to flows. Our main aim in
  the article is to describe the representations in
  Hitchin's component as holonomies of these structures. As
  a preliminary, we show these holonomies form an open set
  in the space of representations.
\item \ref{2}: {\bf Quasi-Fuchsian and Anosov
    representations.} We specify this geometric structure to
  study rank 1 subgroups of semi-simple Lie groups, and more
  specifically the irreducible $PSL(2,\mathbb R)$ in
  $PSL(n,\mathbb R)$; we introduce {\em quasi-Fuchsian
    representations} as deformations of Fuchsian
  representations. In quite a similar way as for classical
  quasi-Fuchsian representations in $PSL(2,\mathbb C)$,
  limit curves appear; instead of taking  values in
  $\mathbb C\mathbb P^{1}$ as in the ``classical`` case,
  they take their values in the flag manifold.  Their
  properties will play a central role in the sequel.
\item \ref{state}: {\bf Statement of the main results.} With
  all the basic notions in hand, we can state our main
  results. The first one is basically that the limit curve
  of a quasi-Fuchsian representation is built from a
  hyperconvex curve. The second one is that every
  representation in Hitchin's component is quasi-Fuchsian.
  We also state converse results.
\item \ref{6}: {\bf Hyperconvex curves.} In this section, we
  study more specifically hyperconvex curves and prove in
  general they admit ``left'' and ``right'' osculating
  flags. This section is independent of the rest of the
  article.
\item \ref{7}: {\bf Preserving a hyperconvex curve.} We
  prove that a representation preserving a hyperconvex curve
  is the holonomy of an Anosov structure.
\item \ref{4}: {\bf Curves and Anosov representations.} This
  is the core of the the article: we prove the limit curve
  of certain Anosov representation is the osculating flag of
  a Frenet hyperconvex curve in Corollary
  \ref{mainlemmacurve3}.
\item \ref{5}: {\bf Anosov representations, 3-hyperconvexity
    and Property (H)} In the core of the proof of the
  previous result, certain properties to be satisfied by
  limit curves were introduced. We study here their
  relations with quasi-Fuchsian representations.
\item \ref{8}: {\bf Closedness.} We show the set of
  quasi-Fuchsian representations is closed in the space of
  all representations. This helps us to conclude the proof
  of our main results.
\item \ref{9}: {\bf Appendix: some lemmas.}
\end{itemize}
\subsection{Further discussions and conjectures}
This section is rather programmatic, it contains some
announcements, precise conjectures. This
section should be skipped by a reader interested in concrete
results.  It is a rather random collection of remarks which
is aimed at suggesting many aspects of Teichmüller theory
considered as a dictionary between various fields of
mathematics should extend to these Hitchin's components.

\subsubsection{Crossratios: $n=\infty$}\label{ninfty}
In a subsequent article, currently under preparation
\cite{FL3}, we explain the relation between Hitchin's
component and crossratios on $\grf$. We define a
crossratio on $\grf$ is a real Hölder function $b$
defined on $(\partial_{\infty}\grf)^{4}\setminus \{(x,y,z,t)/x=w, z=y\}$
satisfying the following rules
\begin{eqnarray*}
b(x,y,z,t)&=&\frac{b(x,y,z,w)}{b(x,w,z,t)},\\
b(x,y,z,t)&=&b(x,t,z,y)^{-1},\\
b(x,y,z,t)&=&b(z,t,x,y).
\end{eqnarray*}
As an example of crossratio, one has the classical
projective crossratio and the crossratio associated by J.-P.
Otal to a negatively curved metric on $S$ \cite{JPO}. They were extensively studied by U. Hamenstädt in \cite{UH} (Notice however our definition includes more general crossratios than those she defined and that some of her results are not true in our generality).  For a
complete description of various aspects of crossratio, one
is advised to read F. Ledrappier's presentation \cite{led}.
Associated to a crossratio are numbers called {\em periods}.
If $\gamma$ is an element of $\pi_1 (S)$, let $\gamma^{+}$
(resp. $\gamma^{-}$) be the attracting (resp.  repelling)
fixed point of $\gamma$ on $\partial_\infty(\pi_1 (S))$. We define the period $l(\gamma)$ of $\gamma$ by
$$
\forall y\in \partial_\infty(\pi_1 (S)),\ \ 
l(\gamma)=\log\vert b(\gamma^+,\gamma^-,y,\gamma y)\vert.
$$
It turns out that a crossratio is completely determined
by its set of  periods which in the case of Otal's crossratio is 
just the collection of lengths of the corresponding closed orbits.

The main result of our article \cite{FL3} explains that
there exists a correspondence between representations in
Hitchin's component and crossratios satisfying some
functional relations, one for each $n$, which are
completely explicit but technical to state. Under this
correspondence, the period of $\gamma$ is equal to
$$
\log(\frac{\lambda_{max}(\rho(\gamma))}{\lambda_{min}(\rho(\gamma))}),
$$
where $\rho$ is the corresponding representation and
$\lambda_{max}(A)$ (resp. $\lambda_{min}(A)$) is the largest
(resp. smallest) real eigenvalue of the matrix $A$.
According to this result, each component
$Rep_{H}(\grf,PSL(n,\mathbb R))$ embeds in the space of all
crossratios, which may be consider as a candidate for
$Rep_{H}(\grf,PSL(\infty,\mathbb R))$. Notice also that the
space of all crossratios is identified ({\em cf.}
\cite{led}) as the space of all Hölder Anosov flows on
the unit tangent bundle of the surface.  This is a rather
mysterious picture, but is has the advantage of (almost)
describing Hitchin's component as a space of objects,
crossratios or Hölder Anosov flows, that may be thought as
``geometric structures'' on the surface.

Talking with Hitchin, we also realised this picture is
coherent with a conjectural picture of his. Namely, he
suggested to consider the group $SL(\infty,\mathbb R)$ as
the group of symplectic diffeomorphisms of $\mathcal G=
\mathbb{RP}^{1}\times \mathbb{RP}^{1}\setminus \Delta$. On
our dynamical side, a crossratio defines a measure,
equivalent to Lebesgue, on $\mathcal G$. For instance, the
choice of a negatively curved metric defines a symplectic
form on the space of geodesics by symplectic reduction, and
this space is identified with $\mathcal G$ via the
identification of $S^{1}$ with the boundary at infinity. Of
course, this measure is invariant under the action of
$\grf$. It follows that after the conjugation by a
homeomorphism sending the measure associated to the
crossratio to the ``standard measure'' on $\mathcal G$, we
obtain a representation of $\grf$ in the group of symplectic
homeomorphisms of $\mathcal G$.  It is striking that these
two pictures coming from different areas of mathematics
agree.

\subsubsection{Universal Hitchin's component: $g=\infty$}
One should notice that Theorem \ref{mainA} allows us to let
the genus $g$ of the surface goes to infinity an thus
provides an extension of the theory of universal
Teichmüller space. Indeed, we may consider the space  $\mathcal
T(n)$ of all Frenet hyperconvex curves in
$\mathbb P(\mathbb R^{n})$, this is a natural candidate for
the {\em universal Hitchin component}, generalising the
group of quasi-symmetric homeomorphisms when $n=2$. Here are some natural quastions : how sit
the various components in this space ? Does it have a
Kähler geometry ? 

\subsubsection{Frenet curves and integrable systems}
We may now hope to relate the subject with integrable
systems. We strongly suggest to read G. Segal very clear
exposition \cite{GS}. A way to build at least locally a
hyperconvex Frenet curve is through differential equations.
Namely, let's consider a $n^{th}$-order linear differential
operator - a {\em Hill operator} - of the following form
\begin{equation}
L(f)=f^{(n)}+a_{2}f^{(n-2)}+a_{3}f^{(n-3)}+ \ldots a_{n}.\label{Hill}
\end{equation}
If $(f_{i},\ldots,f_{n})$ are $n$ independent solutions of
the equation $L(f)=0$, the projective coordinates given by
$$
[f_{1},\ldots,f_{n}]
$$
defines locally a hyperconvex Frenet curve. A different
choice of $f_{i}$ yields the same curve up to a projective
transformation. Since the curves in Theorem \ref{mainA} have
low regularity (they are usually only $C^{1}$), they cannot
be related to smooth regular operators like the one in
Formula (\ref{Hill}). However one would like to know if they
can described by some operator in a weak sense.

The motivation for this question is the following: the space
of Hill's operator is naturally a symplectic manifold and
its Poisson algebra relates to the so-called
$W(n)$-algebras, where $W(2)$ is the Virasoro Algebra ({\em
  cf} \cite{GS})

Apparently, physicists tend to believe a Teichmüller
theory should hold for these $W(n)$-algebras for which
Hitchin's component would play the role of Teichmüller
space. Honestly, I have never understood in the papers that
allude to this question what they really expect as a link
between $W(n)$-algebras and Hitchin's component. Apparently,
the goal is rather to obtain Hitchin's component as a
''double quotient'' of $W(n)$-algebras like it has been done
by M. Kontsevich for Virasoro algebra \cite{MK}, than to
copy the relation of Virasoro algebra with the universal
Teichmüller space as is provided by our previous
discussion.

But at least Theorem \ref{mainA} provides a relation between
$W(n)$-algebras and Hitchin's component which may well be
coherent with the expected picture. Also, the fact that we
still have a candidate to be a companion for $W(\infty)$ as
discusses in the Paragraph \ref{ninfty} seems appealing.

\subsubsection{Holomorphic differentials and the link with Hitchin's theory}
To prove his theorem, N. Hitchin gave an explicit parametrisation of Hitchin's
component. Namely after the choice of a complex structure
$J$ on the compact surface $S$, he identified the component
$Rep_H(\grf,PSL(n, \mathbb R))$ with the vector space
$$
\mathcal Q(2,J)\oplus\ldots\oplus\mathcal Q(n,J),
$$
where $\mathcal Q(p,J)$ denotes the space of holomorphic
$p$-differentials on the Riemann surface $S,J$.  The main
idea in the proof is to identify representations with
harmonic mappings as in K. Corlette's seminal paper
\cite{KC}, or in \cite{D} and \cite{FL1}. Now a harmonic
mapping $f$ with values in a symmetric space gives rise to
holomorphic differentials $q_{2}(f),\ldots $ in quite a
similar fashion as a connexion gives rise to differential forms
in Chern-Weil theory.

Can one improve this parametrisation, and in particular get
rid of the choice of a complex structure and obtain a
parametrisation by holomorphic objects invariant under the
mapping class group ? Here is a suggestion.  A rather
standard check shows that the quadratic differential part
$q_{2}(f)$ vanishes exactly when $f$ is minimal. We may now
wonder if, fixing the representation $\rho$, we can choose
in a unique manner a complex structure on $S$ so that the
associated harmonic is actually minimal.  Another way to
state this question is the following conjecture which I have
discussed many times with Bill Goldman

\begin{conjecture}
  Let $\rho$ be a representation in Hitchin's component. For
  every complex structure $j$ in Teichmüller space
  $\mathcal T$, let $e(j)$ be the energy of the
  corresponding harmonic mapping. Then $e$ has a unique
  minimum.
\end{conjecture}   

This conjecture is well known to be true for $n=2$. For
$n=3$, one can prove it using the ideas linking real
projective structures, affine spheres, Blaschke metric as in
J. Loftin paper \cite{JL} or in the preprint \cite{FL2}; in
order to complete the circle of ideas contained in these
papers, one has just to realise that, for an affine sphere
$S$, the Blashke metric, seen as a map from $S$ to
$SL(3,\mathbb R)/SO(3)$, is minimal.  For a general $n$, on can at least 
show that $e$ is proper \cite{FL4}.

If the last conjecture is true, using our previous
discussion, we would have proved the following result, which
helps to understand the action of the mapping class group
$\mathcal M (S)$ on Hitchin's component
\begin{conjecture}
  The quotient $Rep_{H}(\grf,PSL(n,\mathbb R))/\mathcal M
  (S)$ is homeomorphic to the total space of the vector
  bundle $E$ over Riemann moduli space, whose fibre at a
  point $J$ is
  $$
  E_{J}=\mathcal Q(3,J)\oplus\ldots\oplus\mathcal Q(n,J).
  $$
\end{conjecture}   
Again, by the previous discussion this result is true for
$n=2$ and $n=3$. The fact that the energy is proper would say that the map we can define using Hitchin's identfication (described in the beginning of this paragraph)  from $E$ to $Rep_{H}(\grf,PSL(n,\mathbb R))/\mathcal M
  (S)$ is at least surjective.

\subsubsection{Compactification}

W. Thurston \cite{FLP} gave a compactification of
Teichmüller space, which has been extended in many ways.
More specifically, A. Parreau gave a compactification of the
set of discrete representations in $SL(n,\mathbb R)$
\cite{AP}. Since all representations in Hitchin's component
are discrete, in particular her work gives a
compactification of $Rep_{H}(\grf,PSL(n,\mathbb R))$. It
should be interesting to relate this compactification and
Theorem \ref{mainA}.

\subsubsection{Further extensions and questions}

So far this article only deals with the group $PSL(n,\mathbb
R)$, although Hitchin's Theorem \ref{Hitch2} actually
extends to adjoint groups of all real split forms. It is
rather tempting to conjecture that at least Theorem
\ref{mainB} extends to this general context. It is actually
obvious in cases like $PSO(n,n+1)$ when the corresponding
component is a subset of the component for $PSL(n,\mathbb
R)$.

Another natural extension is  to consider surfaces with
marked points, the holonomy around marked points being
forced to preserve a full flag. To my knowledge, even the
case $n=3$ is not known, although Hitchin's version has been
worked out \cite{IND}.

Notice however that in their remarkable paper \cite{GB}, Volodya Fock and Sacha Goncharov gave a construction and  a combinatorial description of a ``Teichmüller space" for surfaces with punctures or boundary, as well as coordinates  and Poisson structures. Actually their picture extends to the case of real split group. For the moment, it is not clear yet that their Teichmüller space is indeed is connected component. But it is quite believable.

\section{Geometric Anosov Flows}\label{1}
Our starting point is to obtain representations in
$Rep(\grf,PSL(n,\mathbb R))$ as ho\-lo\-no\-mies of
``geometric structures'' associated to flows.  We prove for
these new geometric structures Proposition \ref{Thurslok}
which is an analog in our context of Ehresmann Theorem,
sometimes called Thurston-Lok Holonomy Theorem \cite{LO}\cite{WG}, that states
that the deformation of the holonomy representation for a
compact manifold can be obtained through a deformation of
the structure.

We then explain the example that arises when considering a
rank 1 subgroup of a semisimple group, thus making sense of
a notion of {\em quasi-Fuchsian representation}. We finally
concentrate on the case which is the subject of this paper,
associated to  the irreducible $PSL(2,\mathbb R)$ in
$PSL(n,\mathbb R)$.

As a motivation for our notion of geometric structure, we
begin with a remark. When one defines a $(G,X)$-geometric
structure on a manifold $M$ as an atlas modelled on $X$ with
transition maps in $G$, one requires that $M$ and $X$ have
the same dimension and that the charts are homeomorphisms
(or at least submersions in the case of transverse
structures to foliations), although this is not formally
necessary.  Indeed for instance, if $X$ is allowed to have a
larger dimension, the corresponding ``geometric structure"
would not be rigid enough and too vague to have a useful
meaning. The presence of a flow an a subsequent hyperbolic
hypothesis will allow us to enlarge the definition in this
direction, and still obtain ``rigid'' geometric structures.

Before proceeding to the definition, we recall the
definition of a {\em contracting (or dilating) bundle} over
a dynamical system.

Let $X$ be a topological space equipped with a flow
$\phi_t$. Let $E$ be a topological vector bundle over $X$
such that the action of $\phi_t$ lifts to an action of a
flow $\psi_t$ by bundle automorphisms. Assume $E$ is
equipped with a metric $g$. The bundle $E$ is {\em contracting}
(resp. {\em dilating}), if there exist positive constants
$A$ and $B$, such that for every $u$ in $E$, for every $t$
such that $t > 0$ (resp. $t<0$)
$$
\Vert \psi_t (u)\Vert \leq Ae^{-B\vert t\vert}\Vert u\Vert.
$$

It is useful and classical to remark that if $X$ is compact,
\begin{enumerate}
\item the metric $g$ plays no role;
\item the parametrisation of the flow plays no role either,
  that is if we change the parametrisation of the flow the
  bundle will remain contracting for this new flow.
\end{enumerate}
Therefore to be contracting or dilating over a compact
topological space $X$ is a property of the orbit lamination
$\mathcal L$, the bundle $E$ and the ``parallel transport'' on
$E$ along leaves of $\mathcal L$.
\subsubsection{$(M,G)$-Anosov structure}

Let $M$ be a manifold equipped with a pair of continuous
foliations $\mathcal E^\pm$, whose tangential distributions
are $E^\pm$, and such that
$$
TM=E^+\oplus E^-.
$$
Let $G$ be a Lie group of diffeomorphisms preserving
these foliations.

Let $V$ be a manifold equipped with an Anosov flow $\psi_t$.
Let $\mathcal L$ be the orbit foliation. Let $\tilde V$ be a
Galois covering with covering group $\Gamma$.

We shall say $V$ is {\it $A$-modelled} on $M$ (``A" stands for Anosov), if there exists
a representation $\rho$ of $\Gamma$ in $G$, the {\it
  holonomy representation}, a continuous map $F$ from
$\tilde V$ to $M$, the {\it developing map}, enjoying the
following properties
\begin{itemize}
\item {\it $\Gamma$-equivariance}:
  $$
  \forall \gamma\in\Gamma,~F\circ\gamma=\rho(\gamma)\circ
  F,
  $$
\item {\it Flow invariance}: $$
  F\circ\psi_{t}(x)= F(x), $$
\item {\it Hyperbolicity}: We consider the induced bundle
  $F^\pm =F^*E^\pm$. By the flow invariance, these bundles
  are equipped with a parallel transport along the orbit of
  $\psi_{t}$ (induced for instance by the pull back of any connection on $E^\pm$). By $\Gamma$-equivariance this parallel
  transport is invariant under $\Gamma$. Our last hypothesis
  is that the corresponding lift of the action of $\psi_{t}$
  on $F^+$ (resp. on $F^-$) is contracting (resp.
  dilating).
\end{itemize}
We also say $(V,\mathcal L)$ admits a {\em $(M,G)$-Anosov
  structure}.
\subsubsection{Remarks}
\begin{enumerate}
\item The continuous map $F$ will have in our examples a
  very low regularity. It will only be Hölder.
\item As we shall see in the proof of Proposition
  \ref{Thurslok}, it will turn out that the notion of being
  $A$-modelled is fairly rigid. In other words, if we fix
  the holonomy representation, the only allowed
  infinitesimal transformation of $F$ are translates by
  $\psi_{t}$.
\item One can link this notion to a very classical one. We
  first consider the associated $M$-bundle over $V$ by
  $\rho$, that is $M_{\rho}=(M\times{\tilde V})/\Gamma$
  where the action is the diagonal one. By construction, we
  have now a $\Gamma$-invariant flow $\varphi_{t}$ on
  $M\times{\tilde V}$ given by
  $\varphi_{t}(m,v)=(m,\psi_{t}(v))$.  This flow gives rise
  to a flow $\phi_{t}$ on $N_{\rho}$ lifting $\psi_{t}$.
  Notice now that $F$ gives rise to a flow equivariant
  section of $N_{\rho}$ which we call $\sigma_{F}$. Now our
  hyperbolicity condition just means that $\sigma_{F}(V)$ is
  a hyperbolic subset of $M_{\rho}$ with respect to $\phi_t$.
\end{enumerate}
 
From this last observation and the stability of hyperbolic
sets, we obtain the following Proposition

\begin{proposition}\label{Thurslok}
  Let $M$ be a manifold equipped with a pair of foliations
  as described above.  Let $G$ be the group of
  diffeomorphisms preserving these foliations.  Let $V$ be a
  compact manifold equipped with an Anosov flow $\psi_{t}$.
  Let $\tilde V$ be a Galois covering with covering group
  $\Gamma$.  Let $\cal O$ be the subset of homomorphisms
  $\rho$ from $\Gamma$ to $G$ which are holonomy
  representations of $(M,G)$-Anosov structures.  Then $\cal
  O$ is open.
\end{proposition}
\proof We use the notations of the previous paragraph. We
first have to prove that $\sigma_{F}(V)$ is an isolated
hyperbolic set of $N_{\rho}$. That is, we have to find
an {\em isolating neighbourhood} $U$ characterised by the
property that
$$
\sigma_{F}(V)=\bigcap_{n\in\mathbb Z}\phi^{n}(U).
$$
Remember that $M$ has a local product structure given by
the two foliations $\mathcal E^{\pm}$.

Let's denote by $\pi$ the fibration $M_{\rho}\to V$
described above. We fix for every $x$ in $V$ a complete
Riemannian metric $g_{x}$ on $\pi^{-1}(x)\approx M$
depending continuously on $x$. If $\pi(y)=x$, we consider
$d^{\pm}_{y}$ the associated distance on the leaves
$\mathcal E^{\pm}_{y}$ through $y$ of the foliations
$\mathcal E^{\pm}$. We denote by $B^{\pm}_{y}(\epsilon)$ the
ball of radius $\epsilon$ on $\mathcal E^{\pm}$ centred at
$y$.

Since $M$ has a local product structure, for every $y$, we
can find a real positive number $\epsilon$, such that
\begin{itemize}
\item for every $x$ in $B^{+}_{y}(\epsilon)$, for every $t$
  in $B^{-}_{y}(\epsilon)$, the leaves $\mathcal E^{-}_{x}$
  and $\mathcal E^{+}_{z}$ have a unique intersection
  $G_{y}(x,z)$ in the ball of centre $y$ and radius
  $10\epsilon$,
\item furthermore $G_{y}$ is a differentiable embedding.
\end{itemize}
We set
$$
U_{y}(\epsilon)= G_{y}\big(B^{+}_{y}(\epsilon)\times
B^{-}_{y}(\epsilon)\big).
$$
 
Since $\sigma_{F}(V)$ is compact, we can find $\epsilon$
that satisfies the above conditions for all $y$ in
$\sigma_{F}(V)$. We now consider
$$
U(\epsilon)=\bigcup_{y\in\sigma_{F}(V)}U_{y}(\epsilon).
$$
This last set is a neighbourhood of $\sigma_{F}(V)$.

For $\epsilon$ small enough, since $\sigma_{F}(V)$ is an
hyperbolic set, there exists positive constants $A$ and $B$,
such that we have
$$
\forall z,w \in B^{\pm}_{y}(\epsilon),\forall t >0, \ \ 
d_{y}(\phi_{\pm t}(z),\phi_{\pm t}(w))\leq d_{y}(z,w)Ae^{-B
  t},
$$
this last condition implies that $U$ is an isolating
neighbourhood.

By Theorem 7.4 of C. Robinson's book \cite{CR}, which is
stated for diffeomorphisms but whose proof extends to flows
from the discussion of the next page, we deduce that
$\sigma_{F}(V)$ is stable. In our case, this implies that
after a small perturbation $\hat\rho$ of $\rho$, there
exists a hyperbolic set $W$ of $N_{\hat\rho}$ a
homeomorphism $h$ from $\sigma_{F}(V)$ to $W$ close to the
identity and conjugating the flows up to a small time
change.

Now, we prove there exists a section $\hat\sigma$ such that
$W=\hat\sigma(V)$. Indeed, $H=\pi\circ h\circ\sigma_{F}$ is
a mapping from $V$ to $V$, $C^{0}$-close to the identity and
conjugating the flows up to a small time change. Since the
flow of $\psi_{t}$ on $V$ is Anosov, we deduce that $H$ is
an homeomorphism. It follows that $\pi: W\to V$ is a
homeomorphism. Hence, $W$ is the image of a section
$\hat\sigma$.

Finally, we know that $W$ is a hyperbolic set; recall that
the tangent spaces to the foliations $\mathcal E^{\pm}$ are
invariant by the flow, it follows that these tangent spaces
remain contracting and dilating bundles after a small
perturbation.  \qed

\section{Quasi-Fuchsian and Anosov representations}\label{2}

We now give concrete examples of the situation described
above.

\subsection{Rank 1 subgroups and geometric Anosov structures}

Let $G$ be a semi-simple group and ${\mathcal G}$ its Lie
algebra. Let $H$ be a connected rank 1 semi-simple subgroup
of $G$.  Associated to this situation, we are going to
describe geometric Anosov structures carried by the unit tangent 
bundle of the symmetric space associated to $H$ with its
geodesic flow.

We introduce some notations.
\begin{itemize}
\item Let $A$ be the real split Cartan subgroup of $H$ and
  $Z(A)$ the centraliser of $A$ in $G$. Let $Z_{0}(A)$ be
  the connected component of $Z(A)$ containing the identity.
\item Write $Z_{0}(A)=U$. Let $M=G/U$. Notice that $G$ acts
  on the left on $M$.
\item Let $U\cap H=W\times A$, where the Lie algebra of $W$
  is orthogonal to $A$.
\item Notice that the right action of $A$ on $H/W$ is
  identified the geodesic flow of the unit tangent  bundle of the
  symmetric space of $H$. Let $\mathcal L$ be the orbit
  foliation of this flow.
\item Let $P^+$ (resp. $P^-$) be the parabolic subgroup
  whose Lie algebra is generated by the eigenvectors of non
  negative (resp. nonpositive) eigenvalues of $ad(A)$.
  Notice that $M$ is an open set in $G/P^+\times G/P^-$. Let
  $\mathcal E^\pm$ be the pair of foliations coming from
  this product structure on $M$.
\end{itemize}

We are interested in $(M,G)$-Anosov structures, which we
abusively call again {\em $(H,G)$-Anosov structures}.

\subsubsection{Fuchsian representations.}
 Our
initial result is the following

\begin{proposition}\label{fuchs}
  Let $\Gamma$ be a torsion free discrete subgroup of $H$.
  Let $V=\Gamma\backslash H/W$. Then $(V,\mathcal L)$ admits
  a canonical $(H, G)$-Anosov structure. The developing map
  is the identification of $H/W$ with the left orbit of $H$
  of the identity class in $M$. The corresponding holonomy
  representation is the injection of $\Gamma$ in $G$ through
  $H$. We call such a representation an {\em $(H,G)$-Fuchsian
    representation}.
\end{proposition}

\proof We let $H$ act on the left on $M=G/U$.  Let $m_{0}$
be the class of the identity in $M$.  We define $F$ from
$H/W$ to $M$ by
$$
F(g)=gm_{0}.
$$
We consider the pulled back vector bundle $E$ on $H$
defined by
$$
E=F^{*}TM.
$$
We also consider the bundles $E^\pm$ that come from the
product structure on $M$.  We wish to prove that the right
$A$ action on $E^\pm$ is contracting/dilating.  Notice that
$H$ acts on the left on $E$ by an action that lifts the
standard left action of $H$ on $H/W$. We denote by $g_{*}$
the linear map from $E_{m_{0}}$ to $E_{gm_{0}}$ associated to the action of an element $g$ of $H$.

Recall that $W$ is compact and that $Wm_{0}=m_{0}$. We can
now choose a metric $q_{m_{0}}$ on $E_{m_{0}} $ invariant by
the action of $W$.  We equip now the bundle $E$ with the
metric $q$ defined by
$q_{gm_0}(u,u)=(g^{*}q_{m_0})(u,u)=q_{m_0}(g^{-1}_{*}(u),g^{-1}_{*}(u))$.
This is a well defined metric.

Notice that this metric is invariant by the left action and
hence by the action of any discrete subgroup. Furthermore,
every left $H$-invariant metric on $E$ arises from this
construction.

We finally have a right action of $A$ on $M$ commuting with
the left $H$ action, this action of $A$ preserves globally
the orbit $F(H/W)$; the corresponding action of $A$ on $H/W$
is the geodesic flow, whenever $H/W$ is identified with the
unit tangent bundle of the symmetric space of $H$. We
therefore obtain a right action of $A$ on $E$. If $a$ is an
element of $A$ and $q$ a left $H$ invariant metric, $\tilde
q=a^{*}q$ is also a $H$-invariant metric completely
determined by $q$. By construction, we have $\tilde
q_{m_{0}}=Ad(a)q_{m_{0}}$, it follows the action of $A$ on
$E^\pm$ is contracting/dilating.  \qed 

\subsubsection{Anosov, quasi-Fuchsian representations and limit curves.}
 Assume
now that $\Gamma$ is a cocompact lattice. We define a {\em
  $(H,G)$-Anosov representation} of $\Gamma$ in $G$ as the
holonomy of a $(H,G)$-Anosov structure on $\Gamma\backslash
H/W$ with its geodesic flow.  We define a {\em
  $(H,G)$-quasi-Fuchsian representation} in $G$ as a
representation in the connected component of Fuchsian
representations in the set of $(H,G)$-Anosov representations
in the space of all representations. From Proposition
\ref{Thurslok}, the set of $(H,G)$-Anosov representations is
open. One can check that $(PSL(2,\mathbb R),PSL(2,\mathbb
C))$-quasi-Fuchsian representations coincides with quasi-Fuchsian
representations in the classical sense. Recall that, in this classical case, a
quasi-Fuchsian representation preserves a quasi circle on
${\mathbb CP}^1$. We explain now the counterpart of this
fact.
\begin{proposition}\label{quasimapping}
  Let $\Gamma$ be a cocompact lattice in $H$.  Let $\rho$ be
  a $(H,G)$-Anosov representation of $\Gamma$ in $G$. Let
  $\partial_{\infty}\Gamma$ be the boundary at infinity of
  $\Gamma$. Then, there exist Hölder $\rho$-equivariant
  mappings $\xi^\pm$ from $\partial_\infty\Gamma$ to
  $G/P^\pm$ called the {\em positive and negative limit curves of $\rho$}, such that  
\begin{itemize}
\item if $x\not=y$, $\xi^+(x)$ and $\xi^-(y)$ are opposite
  parabolics,
\item finally, if $\gamma^+$ is an attractive fixed point of
  $\gamma$ in $\partial_\infty\Gamma$, then $\xi^\pm
  (\gamma^+)$ is an attractive fixed point of $\rho(\gamma)$
  in $G/P^\pm$.
\end{itemize}
\end{proposition}

\proof By definition of an Anosov structure, the stable and
unstable manifold of $\phi_t$ (along $F(\tilde V)$) are the
right and left orbit foliations by $P^+$ and $P^-$.
Furthermore, these foliations are well known to be Hölder
(Theorem 19.1.6. of \cite{HK})

We therefore have $\rho$ equivariant Hölder maps from
$\tilde V$ to $G/P^+$ and $G/P^-$ constant along the central
stable (resp. unstable) foliations of the geodesic flow of
$H/W$. Since the space of these central stable leaves is
identified with $\partial_\infty \Gamma$, we have proved
Proposition \ref{quasimapping}. The final two statements are
immediate.

\qed

\subsection{Irreducible $PSL(2,\mathbb R)$ in $PSL(n,\mathbb R)$.}
From now on we concentrate on the following example. First,
we shall consider $V=US$ the unit tangent bundle of a compact
hyperbolic surface, and $\tilde V$ to be the unit tangent bundle
of its universal cover.  We consider the lamination
$\mathcal L$ given by the orbit foliation of the geodesic
flow. We consider also $\mathcal F^\pm$ the central stable
and unstable foliations of the geodesic flow. It is well
known that this data just depends on the fundamental group
$\grf$ of the surface. Indeed we can describe this data the following way. Let
$$
\Delta_3=\{(x_1,x_2,x_3)\in(\partial_{\infty}\grf)^{3}/\exists
i\not=j, \ x_i=x_j\}.
$$
Let's choose an arbitrary orientation on
$\partial_{\infty}\grf$ and let $\partial_{\infty}\grf^{3+}$
be the space of positively ordered triples. Then the following identification holds
\begin{eqnarray*}
\tilde V&=&\partial_\infty\grf^{3+}\setminus\Delta_3\cr
V&=&(\partial_\infty\grf^{3+}\setminus\Delta_3)/\grf.
\end{eqnarray*}
Furthermore, for every $x=(x_+,x_-,x_0)$ in $US$, the leaf $\mathcal L_x$
of $\mathcal L$ through $x$ in $\tilde V$ is
$$
\mathcal L_x=\{(y_+,y_-,y_0)/ y_+=x_+,\\ y_-=x_-\}.
$$
Similarly
$$
\mathcal F^\pm_x=\{y_+,y_-,y_0)/ x_\pm=y_\pm\}.
$$
We are going to model these flows on a specific
situation, namely we consider
\begin{itemize}
\item $G=PSL(n,\mathbb R)$,
\item $H$ the image of the irreducible representation of
  $PSL(2,\mathbb R)$.
\end{itemize} 
In order to simplify our notation we shall speak of {\em
  $n$-quasi-Fuchsian representations} (resp. $n$-Anosov
structures) or just quasi-Fuchsian representations if the
context is clear, instead of $(H,PSL(n,\mathbb R))$-quasi-Fuchsian representations.

\subsubsection{Description of the model}
In this case, $A$ lies in the interior of the Weyl Chamber
and $U=Z_{0}(A)$ is nothing else than the full Cartan
subgroup of $G$, that is the subgroup of diagonal matrices
in a given basis.  It is useful to think of $M=G/U$ as an
open set in $Flag \times Flag$, where $Flag$ is the space of
flags.

Recall that a point of $M$ is a family of $n$ lines $\mathbb
L=\{L_i\}_{i\in\{1,\ldots,n\}}$ in a direct sum.

\subsubsection{A vector bundle description of $n$-Anosov representations}

We immediately have
\begin{proposition}\label{vectorbundledescription}
  Let $\rho$ be a $n$-Anosov representation of $\grf$ in
  $PSL(n,\mathbb R)$ which can be lifted to $SL(n,\mathbb R)$. Let $E$ be the associated $\mathbb
  R^n$ bundle over $V=US$ with its flat connection $\nabla$.
  Then $E$ splits as the sum of $n$ continuous line bundles
  $V_i$ parallel along the leaves of $\mathcal L$.
  Furthermore, let $E^+=(E^+_i)$ (resp. $E^-=(E^-_i)$) be the
  corresponding positive and negative flag bundles,
\begin{eqnarray*}
E_i^+&=&\bigoplus_{j=1}^{j=i}V_j\cr
E_i^-&=&\bigoplus_{j=n-i-1}^{j=n}V_j.
\end{eqnarray*}
Then, $E^+_i$ (resp. $E^-_i$) is parallel along $\mathcal
F^+$ (resp. $\mathcal F^-$).  Finally, if we lift the action
of $\mathcal L$ by the connection, this action is
contracting on $V_i^*\otimes V_j$ for $i>j$.

Furthermore, if we lift the vector bundle $E$ over $\tilde
V$ and identify this bundle with the trivial bundle $\mathbb
R^{n}\times \tilde V$ using the flat connection, we have the
following identification with the positive and negative limit curves of $\rho$
\begin{eqnarray}
E^{\pm}_{(x_{+},x_{0},x_{-})}&=&\xi^{\pm}(x_{\pm}).\label{identif}
\end{eqnarray}

Conversely, the holonomy of such a connection is an
$n$-Anosov representation.
\end{proposition} 

\proof It suffices to remark that $\mathbb
L=(V_{1},\ldots,V_{n})$ is a section of the associated
$M=PSL(n,\mathbb R)/U$ bundle and the tangent spaces to the
associated foliations are
\begin{eqnarray*}
E^{+}&=&\bigoplus_{i>j} ( V_i^*\otimes V_j ),\\
E^{-}&=&\bigoplus_{i<j} ( V_i^*\otimes V_j ).
\end{eqnarray*}
\qed

\subsubsection{Faithfulness and discreteness}

Recall that an element of a semi-simple Lie group is {\em
  purely loxodromic} if it is conjugate to an element in the
interior of the Weyl chamber. In the case of $PSL(n,\mathbb
R)$, this just means that it is real split with 
eigenvalues of multiplicity 1.

Some people may prefer to call purely loxodromic elements {\it strictly hyperbolic}. However, we feel this last terminology may be confusing from the dynamical systems point of view: purely loxodromic element may well have 1 as an eigenvalue and this is a fact which is not felt to be compatible with strict hyperbolicity  for a dynamicist.

We then have

\begin{proposition}\label{lengthorbit}
  Let $\rho$ be a $n$-Anosov representation. Then for each
  $\gamma$ in $\grf$ different from the identity,
  $\rho(\gamma)$ is purely loxodromic. In particular $\rho$
  is faithful. Furthermore if $\rho$ is $n$-quasi-Fuchsian,
  it is irreducible and discrete.
\end{proposition} 

\proof An element is purely loxodromic if it has an
attractive fixed point in the space of flags. Therefore the
first assertion of the proposition follows from Proposition
\ref{quasimapping}. Next, $\rho$ is obviously faithful since
a loxodromic element is not trivial.  Irreducibility follows
from Lemma \ref{corohitchin} and discreteness from Lemma
\ref{discrete} which are both proved in an independent appendix. \qed

\subsubsection{Basic properties of limit curves and  2-hyperconvexity. }
If $\rho$ is a $n$-Anosov representation, according to
Proposition \ref{quasimapping}, we deduce two Hölder
mappings $\xi^+$ and $\xi^-$ from $\partial_{\infty}\grf$ to
the the corresponding homogeneous spaces $G/P^+$ and $G/P^-$
which in our case are both identified with the space of
flags. Since for every attracting point $\gamma^+$ in
$\partial_{\infty}\grf$ of some element $\gamma$ in $\grf$,
$\xi^\pm$ is an attracting point of $\rho(\gamma)$ in the
space of flags and such an attracting point is unique for a
loxodromic element, we conclude that
$\xi^+(\gamma^+)=\xi^-(\gamma^+)$ and hence, by density of
the fixed points, that $\xi^+=\xi^-$.

From now on, we shall therefore write
$$
\xi^\pm=\xi=(\xi^1,\xi^2,\ldots,\xi^{n-1}).
$$
Here, $\xi^i$ takes values in the Grassmannian of
$i$-planes in $E=\mathbb R^{n}$. By definition, we have
$$
\forall x\in \partial_{\infty}\grf,\ \xi^i (x)\subset
\xi^{i+1}(x).
$$
The curve $\xi$ will be called the {\em limit curve} of
$\rho$.  Notice that for $x\not=y$, $\xi(x)$ and $\xi(y)$
are transverse flags since they correspond to opposite
parabolics ({\em cf.} Proposition \ref{quasimapping}).
Hence, we have the following property which we shall call in
short {\em 2-hyperconvexity},
\begin{eqnarray}
\forall x,y\in \partial_{\infty}\grf, x\not=y\implies \xi^{p}(x)\oplus\xi^{n-p}(y)=E.\label{mlc0}
\end{eqnarray}
Certainly, the curve $\xi$ cannot be any curve; it has to
have some properties. For instance, in the $(PSL(2,\mathbb
R),PSL(2,\mathbb C))$ situation this is a quasi-circle. It
also follows from S. Choi and W.  Goldman's work that in
the $(PSL(2,\mathbb R),PSL(3,\mathbb R))$ case, the curve
$\xi^1$ is $C^1$ and bounds a convex set \cite{CG}.

\section{Statement of the main results}\label{state}
We state now our main theorem concerning the properties of
the curve $\xi$, which generalises S.-Y. Choi and W.
Goldman's situation.

\subsection{Quasi-Fuchsian representations, limit curves and Hit\-chin's component}
Our main Theorem is a slight refinement of Theorem
\ref{mainA}
\begin{theorem}\label{maincurve} Let $\rho$ be a  representation in Hitchin's component.  Then $\rho$ is quasi-Fuchsian. Furthermore,
  let $$
  \xi=(\xi^1,\xi^2,\ldots,\xi^{n-1}) $$
  be its limit
  curve. Then $\xi^{1}$ is a hyperconvex Frenet curve, and
  $\xi$ is its osculating flag.  Furthermore, for any triple
  of distinct points $(x,y,z)$ of $\partial_{\infty}\grf$
  the following sum is direct,
\begin{eqnarray}
\big(\xi^{k+1}(y)+\xi^{n-k-2}(x)\big)
+(\xi^{k+1}(z)\cap\xi^{n-k}(x))=E\label{maincurve14}.
\end{eqnarray}
\end{theorem}
We recall that saying $\xi$ is the osculating flag of the
hyperconvex Frenet curve $\xi^{1}$ means the following.
\begin{itemize}
\item Let $(x_1,\ldots,x_p)$ be pairwise distinct points of
  $\partial_{\infty}\grf$. Let $p$ be an integer. Let
  $(n_1,\ldots,n_l)$ be positive integers such that $$
  l=\sum_{i=1}^{i=p}n_i\leq n.  $$
  Then, the following sum
  is direct
\begin{eqnarray}
\sum_{i=1}^{i=l}\xi^{n_i}(x_i)\label{maincurve11}
\end{eqnarray}
\item Furthermore for every $x\in \partial_{\infty}\grf$,
\begin{eqnarray}
\lim_{(y_1,\ldots,y_p)\rightarrow x, y_i
{\hbox{\tiny all distinct}}}
(\bigoplus_{i=1}^{i=p}\xi^{n_i}(y_i))=\xi^{l}(x).\label{maincurve12}
\end{eqnarray}
\end{itemize}
As a consequence $\xi$ is completely determined by $\xi^1$,
and $\xi^1$ is a $C^1$ curve.  Theorem \ref{maincurve}
together with Proposition \ref{lengthorbit} give rise to
Theorem \ref{mainB}. Theorem \ref{maincurve} is proved in
Paragraph \ref{proof1}.
\subsection{Converse results}
It turns out that the curve $\xi^1$ contains all the
information needed to reconstruct our geometry.
\begin{theorem}\label{preservehyper}
  Let $\rho$ be a representation of $\grf$ in $SL(E)$.  Let
  $\xi^1$ be a $\rho$-equivariant continuous map from
  $\partial_{\infty}\grf$ to $\mathbb P (E)$. Assume that
  for all distinct points $(x_1,\ldots,x_n)$, we have the following direct sum
\begin{eqnarray*}
\xi^1(x_1)+\ldots+\xi^1(x_n)=E.
\end{eqnarray*}
Then $\rho$ is a $n$-Anosov representation and $\xi^1$ is
the projection in $\mathbb P(E)$ of the limit curve $\xi$ of
$\rho$. Finally, $\xi^{1}$ is a hyperconvex Frenet curve and
$\xi$ is its osculating flag.
\end{theorem}

This Theorem is proved in Paragraph \ref{proof2}.
Unfortunately, we cannot prove every $n$-Anosov
representation is quasi-Fuchsian. To our
present knowledge, the set of $n$-Anosov representation
could well not be connected. However, Olivier Guichard recently proved the following result, which was conjectured in an earlier version of the present paper and which gives a complete geometric characterisation of  Hitchin's
component \cite{OG}

\begin{theorem}{\sc[Guichard]}
  Let $\rho$ be a representation of $\grf$ in $SL(E)$. Let
  $\xi^1$ be a $\rho$-equivariant continuous map from
  $\partial_{\infty}\grf$ to $\mathbb P (E)$. Assume that
  for all distinct points $(x_1,\ldots,x_n)$, we have the following direct sum
\begin{eqnarray*}
\xi^1(x_1)+\ldots+\xi^1(x_n)=E.
\end{eqnarray*}
Then the representation $\rho$ is in Hitchin's component.
\end{theorem}
These two
theorems provide a converse result to Theorem
\ref{maincurve}.

\section{Hyperconvex curves}\label{6}

\subsection{Definition and notations}
Let $\xi$ be a map from an interval $J$ to $\mathbb P (E)$.
We shall say $\xi$ is {\it hyperconvex}, if for all
$n$-uples of distinct points $(x_1,\ldots,x_n)$ we have
$$
\xi(x_1)+\ldots+ \xi(x_n)=E.
$$
As a notation if $p\leq n$, we write if
$X=(x_1,\ldots,x_p)$ is a $p$-uple of distinct points
$$
\xi^{(p)}(X)=\xi^{1}(x_1)\oplus\ldots\oplus \xi^{1}(x_p).
$$
We also say $X<x$ if for all $i$, $x_i<x$. We write
$X\rightarrow x^+$ as a shorthand for $X\rightarrow x, X<x$.
We have a similar convention for $X\rightarrow x^-$.

We shall actually need a refinement of the notion of hyperconvexity in order  to take in account {\em non continuous} maps. Let $\Omega$  be an orientation on  $E$. Let now $\xi$ be a map (not necessarily continuous) from an interval $J$ to $\mathbb P (E)$. We say that $\xi$ is {\em $*$-hyperconvex} if  the following sum is direct,
$$
\xi^{(p)}(X)=\xi(x_1)+\ldots+ \xi(x_p).
$$
Furthermore, we require  there exists a map $\hat\xi$, the {\em lift of $\xi$}, with values in  $E\setminus\{0\}$ such that the  following holds
\begin{enumerate}
\item for all $y$ in $J$, $\hat\xi(y)\in \xi(y)$,
\item for all $n$-uple of distinct increasing  points $X=(x_{1},\ldots,x_{n})$,  we have
$$
\Omega(\hat\xi(x_1),\ldots, \hat\xi(x_n)) \geq 0.
$$
Notice that this last inequality and the first condition actually implies
\begin{eqnarray}
\Omega(\hat\xi(x_1),\ldots, \hat\xi(x_n)) > 0.\label{*hyp}
\end{eqnarray}
The existence of this ``coherent" lift should be understood in the following way : the map $\xi$, though not being continuous, preserves some ordering.
\end{enumerate}
It is also obvious that such a lift exists whenever $\xi$ is continuous. It follows that every hyperconvex curve defined on a (contractible) interval is in particular $*$-hyperconvex.

\subsection{Left and right osculating flags}
The main result of this section is the following

\begin{lemma}\label{mainhyperconvex}
  Let $\xi$ be an $*$-hyperconvex map from $J$ to $\mathbb
  P(E)$.  Assume that the sequence $\{X_m\}_{m\in\mathbb N}$
  (resp.  $\{Y_m\}_{m\in\mathbb N}$ ) converges to $(x_1,\ldots,x_p)$
  (resp. $(y_1,\ldots,y_{n-p})$) with 
  $$
  x_1\leq \ldots\leq x_p<y_1\leq y_2\ldots \leq y_{n-p}.
  $$
  Assume also that
  $\{\xi^{(p)}(X_m)\}_{m\in\mathbb N}$ (resp.
  $\{\xi^{(n-p)}(Y_m)\}_{m\in\mathbb N}$) converges to $F$
  (resp. $G$).
  
Then
\begin{eqnarray}
F\oplus G=E.\label{mh1}
\end{eqnarray}

Furthermore, for every $p$, with $n\geq p\geq 1$, there exist maps $\xi_+^p$ and
$\xi_-^p$ from $J$ to $Gr(p,E)$ such that
\begin{eqnarray}
\lim_{X\rightarrow x^\pm}\xi^{(p)}(X)&=&\xi^p_\pm(x),\label{mh2}\\
\lim_{(z,y)\stackrel{ z\not=y}{\rightarrow} x^\pm} (\xi(z)\oplus\xi^p_\pm(y)))&=&\xi^{p+1}_\pm(x),\\ \label{mh3}
\xi^{p}_\pm(x)&\subset&\xi^{p+1}_\pm(x).
\end{eqnarray}

Finally, if $\xi^p_+=\xi^p_-$, then both maps are continuous and
\begin{eqnarray}
\lim_{X\rightarrow x}\xi^{(p)}(X)&=&\xi^p_\pm(x),\label{mh41}
\end{eqnarray}
 In particular,  if $\xi^1_+=\xi^1_-$, then both maps are equal to $\xi$ and the latter
 is continuous.  
\end{lemma}

We shall begin by some preliminaries concerning {\em
  increasing maps}, then prove the Lemma.
\subsection{Increasing maps}
We make precise some properties of {\em increasing maps}. Let $p$
be some integer. Let $I$ be an oriented interval. We define
$$
I^{(p)}=\{(x_1,\ldots,x_{p})\in I^{p}/ x_i \leq
x_{i+1}\}.
$$
We define partial orderings on $I^{(p)}$ by
\begin{eqnarray*}
(x_1,\ldots,x_{p})&\leq&(y_1,\ldots,y_{p})\hbox{ iff } \forall i, \ x_i\leq y_i,\\
(x_1,\ldots,x_{p})&<&(y_1,\ldots,y_{p})\hbox{ iff } \forall i, \ x_i< y_i,\\
\end{eqnarray*}
We also define
$$
\hat I^{(p)}=\{(x_1,\ldots,x_{p})\in I^{p}/ x_i <
x_{i+1}\}.
$$

Let now $f$ be a map from $\hat I^{(p)}$ to $\mathbb R$. We
say $f$ is {\em increasing} if for every $(x_1,\ldots,x_p)$
in $\hat I^{(p)}$ and for every $j$,
\begin{eqnarray*}
x_{j-1}<z<y<x_{j+1}&\implies& \\
 f(x_0,\ldots,x_{j-1},z,x_{j+1},\ldots,x_p)&\leq& f(x_0,\ldots,x_{j-1},y,x_{j+1},\ldots,x_p).
\end{eqnarray*}
Notice that this immediately implies
$$
X\leq Y\implies f(X)\leq f(Y).
$$
For an increasing map $f$ and $X$ in $I^{(p)}$, we define
\begin{eqnarray*}
f^+(X)&=&\inf_{Y > X} f(Y) \cr
f^-(X)&=&\sup_{Y < X} f(Y).
\end{eqnarray*}
The next proposition summarises the properties that we shall
need in the sequel. All these properties are immediate.
\begin{proposition}\label{increase}
  Assume $f$ is increasing. Then
\begin{eqnarray}
f^-(X)&=&\lim_{Y\stackrel{Y < X}\rightarrow X} (f(Y))\label{fmol2}\\
f^+(X)&=&\lim_{Y\stackrel{Y > X}\rightarrow X}(f(Y))\label{fmol3}\\
f^-(X)&\leq&f^+(X)\label{mol4}\\
f^+(X)&\leq&f^-(Y), \hbox{\em if }X< Y,\label{fmol5}\\
f^-(X)&\geq&f^+(Y), \hbox{\em if }X> Y.\label{fmol6}
\end{eqnarray}
Finally, assume that $f^+=f^-$ are equal everywhere, then
they are both continuous and
$$
\lim_{Y\rightarrow X}f(Y)=f^+(X).
$$
\end{proposition}

\subsection{$*$-Hyperconvex curves and increasing maps}\label{hyperincrease0}

We begin with the following observation which follows at
once from hyperconvexity.  Let
$(y_1,\ldots,y_{n-1},w_1,w_2)$ be distinct points of the
one-dimensional manifold $J$. Write $Y=(y_1,\ldots,y_{n-1})$
and $W=(w_1,w_2)$. Let $\xi$ be a $*$-hyperconvex curve from $J$
to $\mathbb P (E)$. Then
\begin{eqnarray}
\dim\big(\xi^{(n-1)}(Y)\cap \xi^{(2)}(W)\big)&=&1.\label{hyperinj1}
\end{eqnarray}
Let $I$ an interval contained in  $J\setminus\{w_1,w_2\}$.  We now
consider the map $f_W$
\begin{eqnarray*}
\mapping{\hat I^{(n-1)}}{\mathbb P \big(\xi^{(2)}(W)\big)\setminus\{\xi(w_1)\}}{Y=(y_1,\ldots y_{n-1})}{
\xi^{(n-1)}(Y)\cap\xi^{(2)}(W)
}.
\end{eqnarray*}
Notice this map is well defined thanks to Assertion
(\ref{hyperinj1}) and the fact that
$$
\xi(y_1)\oplus\ldots
\oplus\xi(y_{n-1})\oplus\xi(w_{1})=E.
$$
We prove now
\begin{proposition}\label{hyperincrease} Assume that the lift $\hat\xi$ of $\xi$ is well defined on $J$. Then, 
  for a suitable choice of orientation on $I$ and $\mathbb
  P\big(\xi^{(2}(W)\big)\setminus\{\xi(w_1)\}$, the map
  $f_{W}$ is increasing.
\end{proposition}

\proof 
Let $u_1=\hat\xi(w_1)$. If $y_1<y_2<\ldots<y_{n-1}$ are in $I$, we have by Inequality \ref{*hyp}
$$
\Omega(u_1,\hat\xi(y_1),\ldots,\hat\xi(y_{n-1})) >0.
$$
Finally, we choose the orientation on $\xi^{(2)}(W)$
given by the form
$$
\omega(w,t)=\Omega(\hat\xi(y_1),\ldots,\hat\xi(y_{n-2}),w,t).
$$
 Thanks to  Inequality (\ref{*hyp}), we notice this orientation is independent on the choice
of $(y_1,\ldots,y_{n-2})$ in $I$ provided that
$$
y_1<\ldots <y_{n-2}.
$$
This choice gives an ordering of $\mathbb
P\big(\xi^{(2}(W)\big)\setminus\{\xi(w_1)\}$ in the
following way. For every $L$ in $\mathbb
P\big(\xi^{(2}(W)\big)\setminus\{\xi(w_1)\}$, we choose
$x(L)$ in $L$ such that $\omega(u_1,x(L)) >0$. Finally, we
say $L<L^\prime$ if
$$
\omega(x(L),x(L^\prime))>0.
$$
We can now prove that the map $f_W$ is increasing. Let
$$
Q=\xi (y_1)\oplus\ldots\xi (y_{n-2}).
$$
Let $z$ be such that
$$
y_1<y_2<\ldots<y_{j-1}<z<y_j<\ldots<y_{n-2}.
$$
Let 
$$
L_z=f_{W}(y_1,\ldots,y_{j-1}, z,y_{j},\ldots,y_{n-2}).
$$
Let $\hat x(z)$ in $L_z$ be such that
$$
\hat x(z)=(-1)^{n-j-1}\hat\xi(z)+w(z),\ \ w(z)\in Q.
$$
Then $\omega(u_1,\hat x(z)) >0$. Assume now that
$$
y_1<\ldots <y_{j-1}<z<t<y_{j}<\ldots<y_{n-2}.
$$
Then we have,
\begin{eqnarray*}
\omega(\hat x(z),\hat x(t))&=&\Omega(\hat\xi(y_1),\ldots,\hat\xi(y_{n-2}),\hat x(z),\hat x(t))\\
&=&\Omega(\hat\xi(y_1),\ldots,\hat\xi(y_{j-1}),\hat\xi(z),\hat\xi(t),\hat\xi(y_{j}),\ldots,\hat\xi(y_{n-2}))\\
&>&0.
\end{eqnarray*}
We have just proved
$$
f_{W}(y_1,\ldots,y_{j-1},
z,y_{j},\ldots,y_{n-2})<f_{W}(y_1,\ldots,y_{j-1},
t,y_{j},\ldots,y_{n-2}).
$$
\qed
\subsection{Proof of Lemma \ref{mainhyperconvex}}
\subsubsection{First step: Assertion (\ref{mh2})}
\proof We use the notations of the previous paragraph. Let
$p$ be an integer less than $n$.  Let $x$ be a point in $J$,
$I$ a small neighbourhood of $x$ and
$Z=(y_1,\ldots,y_{n-p-1},w_1,w_2)$ be some tuple of
cyclically oriented points in $J\setminus \{x\}$. Write
$Y=(y_1,\ldots,y_{n-p-1})$, $W=(w_1,w_2)$. According to
Propositions \ref{increase} and \ref{hyperincrease}, we
obtain that there exist maps $F^\pm_{p,Y,W}$ from $I$ to
${\mathbb P}(\xi^{(2)}(W))$ such that
\begin{eqnarray*}
\lim_{X\mapsto x^\pm}
\Big((\xi^{(p)}(X)\oplus \xi^{(n-p-1)}(Y))\cap\xi^{(2)}(W)\Big)=F^\pm_{p,Y,W}(x).
\end{eqnarray*}
Using the fact the choice of $Z$ is arbitrary, we will now
show that there exist maps $\xi^p_\pm$ verifying Assertion
(\ref{mh2}) and characterised by
$$
\big(\xi^p_\pm(x)\oplus
\xi^{(n-p-1)}(Y)\big)\cap\xi^{(2)}(W)=F^\pm_{p,Y,W}(x).
$$
Let's prove this last point in detail.

This in done in two steps. First, let's fix $Y$. Let $U$ be
an interval of $J\grf\setminus I\cup Y$.  We
consider now the subspace
$$
H^\pm_{Y}(x)=\sum_{W\in U^{(2)}}F^\pm_{p,Y,W}(x).
$$
Notice that,
$$
\sum_{W\in U^{(2)}}\Big(\big(\xi^{(p)}(X)\oplus
\xi^{(n-p-1)}(Y)\big)\cap\xi^{(2)}(W)\Big)\subset
\xi^{(p)}(X)\oplus \xi^{(n-p-1)}(Y),
$$
hence,
\begin{eqnarray}
\dim\Big(\sum_{W\in U^{(2)}} \big(\xi^{(p)}(X)\oplus
\xi^{(n-p-1)}(Y)\big)\cap\xi^{(2)}(W)\Big)\leq n-1.\label{lrof11}
\end{eqnarray}
We deduce that $\dim (H^\pm_Y)\leq n-1$. We will now
prove that $\dim (H^\pm_Y)= n-1$ and
\begin{eqnarray}
\lim_{X\mapsto x^\pm}
(\xi^{(p)}(X)\oplus \xi^{(n-p-1)}(Y))=H^\pm_{Y}(x).\label{proohyperincrease1}
\end{eqnarray}
Let $\{X_n\}_{n\in\mathbb N}$ be a subsequence converging to
(let's say) $x^+$, such that
$$
P_n =\xi^{(p)}(X_n)\oplus \xi^{(n-p-1)}(Y),
$$
converges to some hyperplane $H$. By hyperconvexity, we
choose $w_1$ in $U$ such that
$$
H\oplus \xi(w_1)=E.
$$
By hyperconvexity again, we choose $(w_2,\ldots,w_n)$ in $U$
such that
$$
\xi(w_1)\oplus\ldots\xi(w_n)=E.
$$
Let $W_i=(w_1,w_i)$. We notice then that
$$
H=\bigoplus_{i} (H\cap \xi^{(2)}(W_i)).
$$
Since
$$
H\cap \xi^{(2)}(W_i)=F^+_{p,W_i,Y}(x).
$$
It follows that $H\subset H^\pm_Y$, hence $\dim (H^\pm_Y)\geq
n-1$. Combining with Inequality (\ref{lrof11}), we obtain that   $H=H^\pm_Y$, hence Assertion (\ref{proohyperincrease1}).

Our next step follows a similar path. We now
consider an interval $U$ not containing $x$, and set
$$
\xi^p_\pm(x)=\bigcap_{Y\in U^{(n-p-1)}}H^\pm_{p,Y}(x).
$$
Since
$$
\xi^{(p)}(X)\subset \bigcap_{Y\in
  U^{(n-p-1)}}(\xi^{(p)}(X)\oplus \xi^{(n-p-1)}(Y)),
$$
we get,
\begin{eqnarray}
\dim \xi^p_\pm(x)\geq p.\label{lrof12}
\end{eqnarray}
We prove now that $\dim \xi^p_\pm(x)\leq p$ and Assertion
(\ref{mh2}). Let again $\{X_n\}_{n\in\mathbb N}$ be a
subsequence converging to (let's say) $x^+$, such that
$\xi^{(p)}(X_n)\oplus \xi^{(n-p-1)}(Y)$ converges to some
$p$-plane $P$.  By hyperconvexity, we now choose
$(y_1,\ldots,y_{n-p})$ in $U$ such that
$$
P\oplus\xi(y_1)\oplus\ldots\oplus\xi(y_{n-p})=E.
$$
Write $Y_i=(\ldots,y_j,\ldots)_{j\not=i}$, and notice
that
$$
P=\bigcap_{i}(P\oplus \xi^{(n-p-1)}(Y_i)).
$$
In particular, since
$$
P\oplus \xi^{(n-p-1)}(Y_i)=H_{p,Y_i}(x),
$$
we get that $\xi^p_\pm(x)\subset P$, hence
$\xi^p_\pm(x)=P$ thanks to Inequality (\ref{lrof12}), hence Assertion (\ref{mh2}).  \qed
\subsubsection{Second step: completion of the proof of Lemma \ref{mainhyperconvex}}

\proof Assume that $\{X_m\}_{m\in\mathbb N}$ (resp.
$\{Y_m\}_{m\in\mathbb N}$) converges to $(x_1,\ldots,x_p)$
  (resp. $(y_1,\ldots,y_{n-p})$) with 
  $$
  x_1\leq \ldots\leq x_p<y_1\leq y_2\ldots \leq y_{n-p}.
  $$
Assume that $\{\xi^{(p)}(X_m)\}_{m\in\mathbb N}$
(resp.  $\{\xi^{(n-p)}(Y_m)\}_{m\in\mathbb N}$) converges to
$F$ (resp. $G$). We want to show
\begin{eqnarray}
F\oplus G=E.
\end{eqnarray}
Let's assume this is not true. We consider now the smallest
integer $m$ for which there exist integers $p$ and $q$, such
that $p+q= m$, satisfying the following property: there
exist sequences $\{X_m\}_{m\in\mathbb N}$ (resp.
$\{Y_m\}_{m\in\mathbb N}$ ) converging to $(x_1,\ldots,x_p)$
  (resp. $(y_1,\ldots,y_{n-p})$) with 
  $$
  x_1\leq \ldots\leq x_p<y_1\leq y_2\ldots \leq y_{n-p}.
  $$
 such that
\begin{itemize}
\item $\{\xi^{(p)}(X_m)\}_{m\in\mathbb N}$ (resp.
  $\{\xi^{(q)}(Y_m)\}_{m\in\mathbb N}$) converges to $P$
  (resp. $Q$);
\item $P\cap Q\not=\{0\}$.
\end{itemize}
Let $H=P+Q$. Write $X_m=(x^m_1,\ldots,x_p^m)$ and
$Y_m=(y^m_1,\ldots,y^m_q)$, for $m$ large enough with
$x^{m}_i<x^{m}_{i+1}<y^{m}_j<y^{m}_{j+1}$. Let's introduce
$$
X^-_m=(x^m_1,\ldots,x_{p-1}^m),\ \ 
Y^-_m=(y^m_2,\ldots,y^m_q).
$$
We can assume, after extracting a subsequence, that
$\{\xi^{(p-1)}(X^-_m)\}_{m\in\mathbb N}$ and
$\{\xi^{(q-1)}(Y^-_m)\}_{m\in\mathbb N}$ converge
respectively to $P^-$ and $Q^-$. By the minimality of
$m=p+q$, we obtain that
\begin{eqnarray}
P^-\oplus Q=P\oplus Q^-=P+Q=H.\label{proomch10}
\end{eqnarray}
Using hyperconvexity, we now choose
$Z=(z_{1},\ldots,z_{n-p-q})$ points in $U$,
$W=(w_{1},w_{2})$ not in $U$ such that the following sums
are direct
\begin{eqnarray}
    H+ \xi^{(n-p-q)}(Z)+ \xi(w_{1})&=&E\label{proomch11},\\
    P^{-}+ Q^{-}+ \xi^{(n-p-q)}(Z)+ \xi^{(2)}(W)&=&E
    \label{proomch12}.
\end{eqnarray}
Using the notations of Paragraph \ref{hyperincrease0}, we
consider now the family of maps $g_{m}$ defined by
$$
g_{m}(t)=f_{W}(X^{-}_{m},t,Y^{-}_{m},Z).
$$
By Proposition \ref{hyperincrease}, all these maps are
increasing.  From (\ref{proomch10}) and (\ref{proomch11}),
we deduce that
$$
\lim_{m\rightarrow\infty}(g_{m}(x_{p}^{m}))
=\lim_{m\rightarrow\infty}(g_{m}(y_{1}^{m}))
=(H\oplus\xi^{(n-p-q)}(Z))\cap\xi^{(2)}(W):=D.
$$
Recall that
$$
\lim_{m\rightarrow\infty}x_{p}^{m}=x_p,\ \ \ 
\lim_{m\rightarrow\infty}y_{1}^{m}=y_1.
$$
Since all the maps $g_{m}$ are increasing, it follows
that for all $t$ in the interval $I$ joining $x_p$ and $y_1$, we
have
$$
\lim_{m\rightarrow\infty}(g_{m}(t))=D.
$$
On the other hand for all $t$ in $I$, from
(\ref{proomch12}), we have
\begin{eqnarray*}
W_{m}(t)&:=&
\xi^{(p-1)}(X^{-}_{m})\oplus \xi^{(q-1)}(Y^{-}_{m})\oplus\xi (t)\oplus 
\xi^{(n-p-q)}(Z)\\
&=&\xi^{(p-1)}(X^{-}_{m})\oplus 
\xi^{(q-1)}(Y^{-}_{m})\oplus \xi^{(n-p-q)}(Z)\oplus g_{m}(t).
\end{eqnarray*}
It follows that for all $t$ in $I$,
$$
\xi(t)\subset\lim_{n\rightarrow\infty}W_{m}(t)=P^{-}
\oplus Q^{-}\oplus \xi^{(m-n-2)}(Z)\oplus D\subsetneq E.
$$
This last assertion contradicts hyperconvexity, hence
finishes the proof of Assertion (\ref{mh1}).

Assertions (\ref{mh1}) and (\ref{mh2}) imply trivially
Assertion (\ref{mh3}). We finally notice that the final
assertion concerning the case where $\xi^p_+=\xi^p_-$ is a
consequence of the last statement of Proposition
\ref{increase}.\qed

\section{Preserving a hyperconvex curve}\label{7}
We prove the following converse of Proposition \ref{vectorbundledescription}.

\begin{theorem}\label{mainconverse}
  Let $\rho$ be a representation of $\grf$ in $SL(E)$. Let
  $\xi^1$ be a $\rho$-equivariant $*$-hyperconvex map from
  $\partial_{\infty}\grf$ to $\mathbb P(E)$. Assume that
  there exist for all integers $p$ less than $n$,
  $\rho$-equivariant maps $\xi^p_\pm$ from  $\partial_{\infty}\grf$ to $Gr(p,E)$ such
  that 
\begin{eqnarray}
\lim_{y\rightarrow x^\pm}(\xi^1(y)\oplus\xi_\pm^p(x))&=&\xi_\pm^{p+1}(x).\label{limi111}
\end{eqnarray}
We also assume that if
\begin{itemize}
\item $\{x_m\}_{m\in\mathbb N}$ (resp. $\{y_m\}_{m\in\mathbb
    N}$) converges to $x$, (resp.  to $y$), with $x\not=y$, 
\item if $p+q\leq n$ and  $Z=(z_1,\ldots,z_{n-p+q})$ are $n-p-q$ points pairwise distinct and different from $x$ and $y$,   
\item $\{\xi_+^p(x_m)\}_{m\in\mathbb N}$ (resp.
  $\{\xi_-^p(x_m)\}_{m\in\mathbb N}$,
  $\{\xi_+^{q}(y_m)\}_{m\in\mathbb N}$) converges to $P^+$
  (resp.  to $P^-$, $Q$)
\end{itemize}
Then
\begin{eqnarray}
P^\pm\oplus Q\oplus \xi^1(z_1)\oplus\ldots\oplus \xi^1(z_{n-p-q})&=&E,\label{sumdir}
\end{eqnarray}
As a conclusion, then
\begin{itemize}
\item $\rho$ is $n$-Anosov,
\item $\xi_{+}^{p}=\xi_{-}^{p}$,
\item $(\xi^1,\xi^2_-,\ldots,\xi^{n-1}_-)$ is the limit
  curve of $\rho$.
\end{itemize}
\end{theorem}
The following corollary is immediate

\begin{coro}
  Let $\rho$ be a representation of $\grf$ in $SL(E)$. Let
  $$
  \xi=(\xi^1,\ldots,\xi^{n-1})
  $$
  be a $\rho$-equivariant continuous map from
  $\partial_{\infty}\grf$ to $Flag(E)$.  Assume $\xi^{1}$ is
  hyperconvex and that
\begin{eqnarray*}
\forall x,y\in\partial_\infty\grf,\forall p,~x\not=y&\implies&\xi^p(x)\oplus\xi^{n-p}(y)=E,\\
\lim_{y\rightarrow x}(\xi^1(y)\oplus\xi^p(x))&=&\xi^{p+1}(x).
\end{eqnarray*}
Then $\rho$ is $n$-Anosov and $\xi$ is the limit curve of
$\rho$.
\end{coro}
As a corollary proved in Paragraph \ref{proof2}, we obtain
Theorem \ref{preservehyper}.
\subsection{Proof of Theorem \ref{mainconverse}}
\proof Let's choose an orientation on
$\partial_{\infty}\grf$. Let
$$
M=\{(x,y,w)\in \partial_{\infty}\grf^3, \hbox{ distinct
  and cyclically oriented }\}.
$$
Notice that $\grf$ acts properly on $M$, in such a way
the quotient is compact and homeomorphic to the unit tangent  bundle
of the surface $S$.  We write the generic element $x$ of $M$
as $x=(x_+,x_0,x_-)$. Consider on $M$ the lamination
whose leaves are
$$
\mathcal L_{x_+,x_-}=\{(x_+,w,x_-)\in M/ w\in
\partial_{\infty}\grf \}.
$$
This lamination is $\grf$ equivariant and its quotient is
identified with the lamination by leaves of the geodesic
flow. Consider also the following 2-dimensional
laminations whose leaves are:
\begin{eqnarray*}
\mathcal F^+_{x_+}&=&\{(x_+,w,y)\in M/ w,y\in \partial_{\infty}\grf \}\\
\mathcal F^-_{x_-}&=&\{(y,w,x_-)\in M/w,y\in \partial_{\infty}\grf\}.
\end{eqnarray*}
Now consider the $E$-associated bundle on $M/\grf$ to
$\rho$, also denoted abusively by $E$. Consider the subbundles
$E^+_i$ and $E^-_i$ of $E$ given by
\begin{eqnarray*}
{E_i^+}_{(x_+,x_0,x_-)}&=&\xi_+^i(x_+)\\
{E_i^-}_{(x_+,x_0,x_-)}&=&\xi_-^i(x_-).
\end{eqnarray*}
Notice these bundles are not {\it a priori} continuous. The
bundle $E_i^+$ (resp. $E_i^-$) are parallel along the leaves of
$\mathcal F^+$ (resp. $\mathcal F^-$).  Let $V^i=E_i^+\cap
E_{n-i+1}^-$. It is a well defined 1-dimensional subbundle
of $E$ (thanks to Hypothesis (\ref{sumdir})), which is parallel
along the leaves of $\mathcal L$. Notice that a subspace
supplementary to $V^i$ is $E^+_{i-1}\oplus E^-_{n-i}$. Denote by $ \alpha_i$ a 1-form whose kernel is that
supplementary subspace.

We first define a metric on $((V^i)^*\otimes V^{i+1})_{w}$.
Write $w=(x_+,x_0,x_-)$. Let $u$ be a nonzero element of
$V^i$, $z(w)$ a a nonzero element of $\xi^1(x_0)$. The
metric on $((V^i)^*\otimes V^{i+1})_{w}$ is given by
$$
\Vert \phi\Vert_w =\Big\vert\frac {\langle \alpha_{i+1}\vert\phi
  (u)\rangle \langle\alpha_{i}\vert z(w)\rangle}{ \langle
  \alpha_{i+1}\vert z(w)\rangle \langle\alpha_{i}\vert
  u\rangle}\Big\vert.
$$
This metric is well defined since by Hypothesis (\ref{sumdir}) the following sum is direct
\begin{eqnarray*}
\xi_+^{j-1}(x_+)+\xi_-^{n-j}(x_-)+ \xi^1(x_0),
\end{eqnarray*}
and in particular for all $j$, 
$$
\langle\alpha_j\vert z(w)\rangle\not=0.
$$
Obviously, this metric is independent of the choice of
$u$ and $z(w)$.  For the moment this metric is not obviously
continuous (although with little effort, one could prove it
is bounded).

We are going to prove now that $(V_i)^*\otimes V_{i+1}$ are
{\em weakly contracting bundles} for this metric.  By weakly
contracting bundle, we mean the following: if $\sigma$ is a
parallel section of $(V_i)^*\otimes V_{i+1}$ along a leaf of
$\mathcal L$, then
\begin{eqnarray*}
\lim_{x_0\rightarrow x_+}\Vert \sigma\Vert_{(x_+,x_0,x_-)} &=&0\\
\lim_{x_0\rightarrow x_-}\Vert\sigma\Vert_{(x_+,x_0,x_-)}&=&\infty.
\end{eqnarray*}
Notice that a parallel section $\sigma$ of $(V_i)^*\otimes
V_{i+1}$ along a leaf of $\mathcal L$, corresponds to a
fixed element $\phi$ in $(V_i)^*\otimes V_{i+1}$.  Then
$$
\frac{\Vert \phi\Vert_w}{\Vert \phi\Vert_t}=\Big\vert\frac
{\langle\alpha_{i+1}\vert z(t) \rangle \langle
  \alpha_{i}\vert z(w) \rangle}{ \langle\alpha_{i+1}\vert
  z(w)\rangle\langle\alpha_{i}\vert z( t)\rangle }\Big\vert.
$$
Let $w=(x_{+},x_{0},x_{-})$ and imagine now that $x_0$
converges to $x_+$.  By Hypothesis (\ref{limi111}), we now
may choose $z^\prime (w) = \alpha +\beta $ in
$$
\xi^1(x_0)\oplus E^{+}_{i-1}=\xi^1(x_0)\oplus
\xi^{i-1}_+(x_+),
$$
with $0\not=\alpha \in\xi^1(x_0)$ and
$\beta\in\xi^{i-1}_+(x_+)$ such that $z^\prime(w)$ converges
to a nonzero vector $u \in
\xi_{+}^{i}(x_{+})\cap\xi_{-}^{n-i+1}(x_{-})$ when $x_0$
converges to $x_+$. Notice that
$$
\frac {\langle\alpha_{i+1}\vert z(t) \rangle \langle
  \alpha_{i}\vert z(w)\rangle}{\langle\alpha_{i}\vert
  z(t)\rangle \langle \alpha_{i+1}\vert z(w)\rangle}=\frac
{\langle\alpha_{i+1}\vert z(t) \rangle \langle
  \alpha_{i}\vert z^\prime(w)\rangle}{\langle\alpha_{i}\vert
  z(t)\rangle \langle \alpha_{i+1}\vert z^\prime(w)\rangle}.
$$
To conclude, we remark that
\begin{eqnarray}
\lim_{x_0 \rightarrow x_+}\frac { \langle \alpha_{i}\vert z^\prime(w)\rangle}{\langle \alpha_{i+1}\vert z^\prime(w)\rangle}
=\frac { \langle \alpha_{i}\vert u\rangle}{\langle \alpha_{i+1}\vert u\rangle}=\infty.
\end{eqnarray}
A similar reasoning when $x_0$ tends to $x_-$, implies the
bundles are weakly contracting.

We can now show that $\xi^p_+=\xi^p_-$. First, by Hypothesis
(\ref{sumdir}), $\xi^p_-(x_{+})$ is the graph of a
homomorphism $\psi$ in $Hom(E^p_+,E^p_-)=(E^p_+)^*\otimes E^p_-$.  Since
$(E^p_+)^*\otimes E^p_-$ is a weakly contracting bundle and
$\psi$ is parallel, to prove $\xi^p_+=\xi^p_-$, it suffices
to show $\Vert\psi\Vert$ is uniformly bounded. Let's prove
that.  Assume we have a sequence of points
$\{x_m\}_{m\in\mathbb N}$ in $M$ such that
$\{\Vert\psi\Vert_{x_m}\}_{m\in\mathbb N}$ tends to
$+\infty$. Since $\grf$ acts cocompactly on $M$ and the
whole situation is invariant under $\grf$, we can as well
assume that $\{x_m\}_{m\in\mathbb N}$ converges to $y$ in
$M$. We can extract a subsequence such that
\begin{itemize}
\item for all $i$, $\{(V_i)_{x_m}\}_{m\in\mathbb N}$
  converges to $W_i$ in $E_y$.
\item $\{(E^i_\pm)_{x_m}\}_{m\in\mathbb N}$ converges to
  $F^i_\pm$,
\item $\{\xi^p_-(x_m)\}_{m\in\mathbb N}$ converges to $Q$ .
\end{itemize}
By Hypothesis (\ref{sumdir}), $Q$ is a graph of a map $\phi$
from $F^p_+$ to $F^p_-$. In particular, $\Vert \phi\Vert$ is
bounded. The contradiction now follows from
$$
\lim_{m\rightarrow\infty}(\Vert
\psi\Vert_{x_m})=\Vert\phi\Vert\not=\infty.
$$

Now that we know $\xi^p_+=\xi^p_-$, we can conclude that
both maps are continuous by Proposition
\ref{mainhyperconvex}. The metric we have defined previously
is therefore continuous. We can repeat the argument above
about the weakly contracting property, the limits that we
obtain  are now uniform, and therefore the argument  shows
the bundles are contracting. Let's do it in more detail. To
prove the bundles $(V_i)^*\otimes V_{i+1}$ are contracting,
we need to show that there exist constant $t_{0}>0$, such
that if $\Psi_{t}$ is the lift of the geodesic flow
$\psi_{t}$ on $M$, then for every vector $\sigma$ in
$(V_i)^*\otimes V_{i+1}$, we have
$$
\forall t >t_{0},\ \Vert \Psi_{t}(\sigma)\Vert\leq
\frac{1}{2} \Vert \sigma \vert.
$$
Let's prove it by contradiction. If this is not true, we
would have a sequence of points $\{w_{n}\}_{n\in\mathbb N}$
of $M/\grf$ a sequence $\{t_{n}\}_{n\in\mathbb N}$
converging to $+\infty$, a sequence
$\{\sigma_{n}\}_{n\in\mathbb N}$ in $((V_i)^*\otimes
V_{i+1})_{w_{n}}$ such that
\begin{eqnarray}
\frac{ \Vert \Psi_{t_{n}}(\sigma_{n})\Vert}{\Vert \sigma_{n} \Vert}\geq \frac{1}{2}.\label{pres11} 
\end{eqnarray}
We may now lift the sequence $w_{n}$ to $M$ and assume,
since the action of $\grf$ is cocompact, that the sequence
converges to $w_{0}=(x^{0}_{+},x^{0}_{0},x^{0}_{-}).$ Let's
write
\begin{eqnarray*}
w_{n}&=&(x^{n}_{+},x^{n}_{0},x^{n}_{-}),\\
\phi_{t_{n}}(w_{n})&=&(x^{n}_{+},y^{n}_{0},x^{n}_{-}).
\end{eqnarray*}
By assumption
\begin{eqnarray*}
\lim_{n\rightarrow\infty}(x^{n}_{+})&=&\lim_{n\rightarrow\infty}(y^{n}_{0})=x^{0}_{+},\\
\lim_{n\rightarrow\infty}(x^{n}_{0})&=&x^{0}_{0},\\
\lim_{n\rightarrow\infty}(x^{n}_{-})&=&x^{0}_{-}
\end{eqnarray*}
Let's denote by $z(t)$ a nonzero element in $\xi^{1}(t)$,
we then have
\begin{eqnarray}
\Big\vert\frac {\langle\alpha_{i+1}\vert z(x_{0}^{n}) \rangle \langle \alpha_{i}\vert z(y_{0}^{n}) \rangle}{ \langle\alpha_{i+1}\vert z(y_{0}^{n})\rangle\langle\alpha_{i}\vert z( x_{0}^{n})\rangle }\Big\vert
=\frac{ \Vert \Psi_{t_{n}}(\sigma_{n})\Vert}{\Vert \sigma_{n} \Vert}\geq \frac{1}{2}.\label{pres12} 
\end{eqnarray}
We now obtain the contradiction knowing that
$$
\lim_{n\rightarrow\infty}(\xi^{1}(y_{0}^{n})\oplus\xi^{i}(x_{+}^{n}))=\xi^{i+1}(x^{0}_{+}),
$$
which implies as above
$$
\lim_{n\rightarrow\infty}(\frac {\langle\alpha_{i+1}\vert
  z(x_{0}^{n}) \rangle \langle \alpha_{i}\vert z(y_{0}^{n})
  \rangle}{ \langle\alpha_{i+1}\vert
  z(y_{0}^{n})\rangle\langle\alpha_{i}\vert z(
  x_{0}^{n})\rangle })=0.
$$
From this it follows the bundle are indeed contracting.
The conclusion finally follows from Proposition
\ref{vectorbundledescription}.  \qed

\subsection{Proof of Theorem \ref{preservehyper}}\label{proof2}

We begin with a lemma of independent interest that will be
used in the sequel, then proceed to the proof

\subsubsection{Direct sums and limits}
Our first lemma is the following:
\begin{lemma}\label{drslsource}
  Let $\xi$ be the limit curve of an Anosov representation.
  Assume, that for all distinct points $(x_1,\ldots,x_q)$ in
  $\partial_{\infty}\grf$, and integers
  $(n_{1},\ldots,n_{q})$ with $k=\sum n_{i}\leq n$, the
  following sum is direct $$
  \xi^{n_{1}}(x_1)+\ldots\ldots\xi^{n_{q}}(x_q); $$
  and
  furthermore
\begin{eqnarray}
\lim_{(x_0,x_1,\ldots,x_l)\rightarrow x}
\big(\xi^{n_{1}}(x_1)
\oplus\ldots\ldots\xi^{n_{q}}(x_q)\big)&=&\xi^{k}(x).\label{drsl1}
\end{eqnarray}
Then, for all $y$ distinct from $(x_1,\ldots,x_q)$, the
following sum is direct
$$
\xi^{n-k}(y)+
\xi^{n_{1}}(x_1)+\ldots\ldots\xi^{n_{q}}(x_q).
$$
\end{lemma}
\proof If $(y,x_{1},\ldots,x_{q})$ is a collection of $q+1$
distinct points of $\partial_{\infty}\grf$, it is a classical fact that  there exist two
distinct points $t$ and $z$, a sequence
$\{\gamma_n\}_{n\in\mathbb N}$ of elements of $\grf$, such
that
\begin{eqnarray*}
\forall i\leq q,\ \lim_n(\gamma_n(x_i))&=&t,\cr
\lim_n(\gamma_n(y))&=&z.
\end{eqnarray*}
Now by 2-hyperconvexity,
$$
\xi^{k}(z)\oplus\xi^{n-k}(t)=E,
$$
for $n$ sufficiently large, by Hypothesis (\ref{drsl1}),
we have
$$
\xi^{n-k}(\gamma_n(y))\oplus\xi^{n_{1}}(\gamma_n(x_1))
\oplus\ldots\oplus \xi^{n_{q}}(\gamma_n(x_q))=E.
$$
Hence the result since
$$
\xi^s(\gamma(w))=\rho(\gamma)\xi^s(w).
$$
\qed
\subsubsection{Proof of Theorem \ref{preservehyper}}

From Lemma \ref{mainhyperconvex} and Theorem
\ref{mainconverse}, we deduce immediately that $\rho$ is
$n$-Anosov. Furthermore, $\xi^{1}$ is the projection in
$\mathbb P(E)$ of the limit curve $\xi$ of $\rho$, and we
have, for every $x\in \partial_{\infty}\grf$,
\begin{eqnarray}
\lim_{(y_1,\ldots,y_l)\rightarrow x, y_i
{\hbox{\tiny all distinct}}}
(\bigoplus_{i=1}^{i=l}\xi^{1}(y_i))=\xi^{l}(x).\label{fre1}
\end{eqnarray}
To conclude the proof of Theorem \ref{preservehyper}, it
suffices to show Assertions (\ref{fre2}) and (\ref{fre3}) of
the definition of Frenet hyperconvex curve.

Let's first prove Assertion (\ref{fre2}).  Let
$(x_1,\ldots,x_p)$ be pairwise distinct points of
$\partial_{\infty}\grf$.  Let $p$ be an integer.  Let
$(n_1,\ldots,n_p)$ be positive integers such that
$$
k=\sum_{i=1}^{i=p}n_i\leq n.
$$
We want to prove the following sum is direct
\begin{eqnarray}
\xi^{n_i}(x_i)+\ldots+\xi^{n_{p}}(x_{p}).\label{fre22}
\end{eqnarray}
We prove it by induction on $p$.  It is true for $p=2$, by
2-hyperconvexity.  Assume it is true for $p=q-1$. Using
Assertion (\ref{fre1}), we deduce that
\begin{eqnarray}
    \lim_{(x_{1},\ldots, x_{q-1})\rightarrow x}(\xi^{n_{1}}(x_{1})
    \oplus\ldots\oplus\xi^{n_{q-1}}(x_{q-1}))=\xi^{k-n_{q}}(x).
    \label{proomaincurve2}
    \end{eqnarray}
    Now Lemma \ref{drslsource} finishes the induction.
    Finally, Assertion (\ref{fre3}) is an immediate
    consequence of Assertions (\ref{fre1}) and (\ref{fre2}).
    \qed
\section{Curves and Anosov representations.}\label{4}

\subsection{Definitions}
We introduce some definitions.
\subsubsection{$(p,l)$-direct.}
Let $(p,l)$ some integers such that $p+l\leq n$. We say the
limit curve $\xi$ (or the corresponding representation
$\rho$) is {\em $(p,l)$-direct} if for all distinct
$(y,x_0,\ldots,x_l)$ the following sum is direct
$$
\xi^{n-p-l}(y)+\xi^p(x_0)+\xi^1(x_1)+\ldots+
\xi^1(x_l).
$$
Notice that to say the representation is $(1,n-1)$-direct
is to say $\xi^{1}$ is hyperconvex.
\subsubsection{$(p,l)$-convergent.}
Let $(p,l)$ some integers such that $p+l\leq n$. We say the
limit curve $\xi$ (or the corresponding representation
$\rho$) is {\em $(p,l)$-convergent} if for all distinct
points $(x_0,\ldots,x_l)$ in $\partial_{\infty}\grf$ the
following sum is direct
$$
\xi^p(x_0)+\xi^1(x_1)+\ldots+ \xi^1(x_l);
$$
and if furthermore
$$
\lim_{(x_0,x_1,\ldots,x_l)\rightarrow
  x}\big(\xi^p(x_0)\oplus\xi^1(x_1)\oplus\ldots\oplus
\xi^1(x_l)\big)=\xi^{p+l}(x).
$$
\subsubsection{3-hyperconvexity}
We say the limit curve (or the corresponding representation)
is {\em 3-hyperconvex}, if $k+p+l\leq n$ and $(x,y,z)$ are
distinct, then the following sum is direct
$$
\xi^{k}(x)+\xi^{p}(y)+\xi^{l}(z).
$$

\subsubsection{Property (H)}
We say the limit curve $\xi$ (or the corresponding
representation) satisfies {\em Property (H)}, if for every
triple of distinct points $x$, $y$ and $z$ and integer $k$,
we have
\begin{eqnarray*}
\xi^{k+1}(y)\oplus(\xi^{k+1}(z)\cap\xi^{n-k}(x))\oplus\xi^{n-k-2}(x)&=&E.
\end{eqnarray*}

\subsection{Main Lemma}
As a main step in the proof of Theorem \ref{maincurve}, we
shall prove the following lemma.

\begin{lemma}\label{mainlemmacurve2}
  Let $\xi$ be the limit curve of an Anosov representation.
  Assume
  \begin{itemize}
  \item it is 3-hyperconvex,
  \item it satisfies Property (H).
  \end{itemize} 
  Then, the curve is $(k,l)$-convergent for all integers
  with $k+l\leq n$. In particular, $\xi^{1}$ is hyperconvex.
\end{lemma}
We state a corollary that follows at once from Theorem
\ref{preservehyper} and which will be the way we shall use
Lemma \ref{mainlemmacurve2}
\begin{coro}\label{mainlemmacurve3}
  Let $\xi$ be the limit curve of an Anosov representation.
  Assume
  \begin{itemize}
  \item it is 3-hyperconvex,
  \item and satisfies Property (H).
  \end{itemize} 
  Then $\xi^{1}$ is a hyperconvex Frenet curve, and $\xi$ is
  its osculating flag.
\end{coro}
We begin with an observation that is a consequence of Lemma
\ref{drslsource}:
\begin{lemma}\label{drsl}
  Let $\xi$ be the limit curve of an Anosov representation.
  Assume the representation is $(p,l)$-convergent, then it
  is $(p,l)$-direct
\end{lemma}

\subsection{Bundles}\label{3}
To prove our Lemma \ref{mainlemmacurve2}, we shall need at
some point the vector bundle description of Proposition
\ref{vectorbundledescription} and prove some preliminary
results in a general situation.
\subsubsection{Hyperconvex rank 2 vector bundle}
Here is the situation we wish to describe.  First let's start with some convention. Let $M$ be a manifold. For any vector bundle $F$ over $M$, we shall denote
by $F_x$ the fibre at a point $x$ of a vector bundle $F$, for any foliation  $\mathcal L$
we denote by $\mathcal L_x$ the leaf 
passing through $x$.

We are first interested in actions of flow on a compact
manifold which preserves a one dimensional foliation. Namely
\begin{itemize}\label{hyp-reg}
\item Let $\phi_{t}$ be a flow of homeomorphisms of a
  compact topological manifold $M$.
\item Let $\mathcal F$ be a 1-dimensional lamination of $M$ with no compact leaves.
  We assume this foliation is invariant under the flow of
  $\phi_t$.
\end{itemize}

This is for instance satisfied when $\phi_t$ is an Anosov
flow and $\mathcal F^+$ is the stable (or unstable)
foliation.

We are interested on bundles over $M$ and actions on these
bundles which lift the previous one. We shall say a vector bundle $E$ over $M$
admits a {\it flag action} if it satisfies the following
assumptions:
\begin{itemize}
\item $E$ is a vector bundle of rank $2$, equipped with a
  parallel transport along leaves of $\mathcal F$;
\item the action of $\phi_{t}$ lift to an action $\psi_{t}$
  by bundle automorphism on $E$, and this action preserves
  the parallel transport;
\item $E$ admits a direct sum decomposition, invariant under $\psi_t$,  into continuous
  oriented subbundles of rank 1
  $$
  E=W^1\oplus W^2, $$
\item $W^2$ is parallel along leaves of $\mathcal F$.
\item {\it Contraction assumption}. We assume $ (W^1)^*\otimes W^2$ is a contracting vector bundle for $\psi_t$. It follows that if we take a 1-dimensional
  vector space $L$ in $W^2_x\oplus W^1_x$ different than
  $W^{2}_{x}$, then
  $$
  \lim_{t\rightarrow\infty}d(\psi_t
  (L),W^1_{\phi_t(x)})=0.
  $$
\end{itemize}

We now introduce some notations.  Assume $x$ and $y$ are on
the same leaf of $\mathcal F$. We shall denote by
$W^1_{x,y}$ the vector subspace of $F^i_x$ which is the
parallel transport of $W^{1}_y$ along the leaf.  We say the
bundle $E$ is {\em hyperconvex} if and only if, for all
distinct $z$ and $y$ in the leaf $\mathcal F_x$
\begin{eqnarray}
W^1_{x,z}\oplus W^1_{x,y}&=&E_x.\label{hyperconvex}
\end{eqnarray}
Our main Lemma is the following
\begin{lemma}\label{rank2}
  Assume the rank 2 vector bundle $E$ equipped with a flag
  action is hyperconvex. Then, the map $J_x$
\begin{displaymath}
\left\{
\begin{array}{rcl}
{\mathcal F}_x&\rightarrow &{\mathbb P}(E_x)\setminus \{W^{2}_x\}\\
y&\mapsto&W^1_{x,y}
\end{array}
\right.
\end{displaymath}
is a homeomorphism; furthermore, for every $x$ in $M$
$$
\lim_{y \rightarrow \infty,\ y\in{\mathcal
    F}_x}(W^1_{xy})=W^{2}_x.
$$
\end{lemma}

We begin by explaining the relations of this notion with our
situation, then prove a preliminary lemma and finally
conclude.
\subsubsection{Hyperconvex bundles and Anosov representations}\label{hbar}
Rank 2 vector bundles with a flag action arise naturally
from Anosov representations. Indeed using the notations of
Proposition \ref{vectorbundledescription},  we shall soon show the bundles
$E^+_k/E^+_{k-2}$ are of this type. More precisely, let $E$
be the vector bundle associated to a $n$-Anosov
representation by Proposition \ref{vectorbundledescription} with its flat connection.
Let $F_{k}=E^{+}_{k}/E^{+}_{k-2}$. Notice that $F$ is
equipped with a flat connection along $\mathcal F^{+}$.
Obviously for this connection
$W^{2}=E^{+}_{k-1}/E^{+}_{k-2}$ is parallel.  We can
identify $F_{k}$ with $E^{-}_{n-k+2}\cap E ^{+}_k$. In this
interpretation, we have
$$
W^{2}= E^{-}_{n-k+2}\cap E^{+}_{k-1}=V^{k-1}.
$$
Let
$$
W^{1}=E^{-}_{n-k+1}\cap E^{+}_{k}=V^{k}.
$$
Then we have the following
\begin{proposition}\label{rkhc}
  Let $E$ be the vector bundle associated to a $n$-Anosov
  representation by Proposition
  \ref{vectorbundledescription}. Then $F_{k}$ with the
  structure described above is a rank 2 vector bundle
  equipped with a flag action. Furthermore, using 
  Identification (\ref{identif}) of Proposition
  \ref{vectorbundledescription}, we have
\begin{eqnarray}
W^{1}_{(x,x_{0},w), (x,x_{0},y)}=
(\xi^{n-k+1}(y) \oplus \xi^{k-2}(x))\cap \xi^{n-k+2}(w)\cap \xi^{k}(x)
\label{identif2}.
\end{eqnarray}
Finally, 
the representation satisfies Property (H) if and only if  the bundles $F_{k}$ are hyperconvex.
\end{proposition}
\proof By definition, the bundle $F_{k}$ (with its flat
connection) along the leaf of $\mathcal F^{+}$ passing
through $(x,x_{0},w)$ is identified with the trivial bundle
whose fibre is
$$
\xi^{k}(x)/\xi^{k-2}(x).
$$
We identify $\xi^{k}(x)\cap\xi^{n-k+2}(w)$ with this fibre using 
the projection along
$\xi^{k-2}(x)$. Then we get
\begin{eqnarray*}
W^{1}_{(x,x_{0},w), (x,x_{0},y)}=
((\xi^{n-k+1}(y)\cap \xi^{k}(x)) \oplus \xi^{k-2}(x))\cap \xi^{n-k+2}(w).
\end{eqnarray*}
This in turn implies  Identification (\ref{identif2}).
Therefore, hyperconvexity is equivalent to, for $y\not=t$
\begin{eqnarray}\label{hyperHproof}
((\xi^{n-k+1}(y)\cap \xi^{k}(x)) \oplus\xi^{k-2}(x))&\cap& \xi^{n-k+2}(w)\cr
\oplus((\xi^{n-k+1}(t)\cap \xi^{k}(x)) \oplus \xi^{k-2}(x))&\cap &\xi^{n-k+2}(w)\cr
=\xi^{k}(x) &\cap& \xi^{n-k+2}(w).
\end{eqnarray}
If we add $\xi^{k-2}(x)$ to both sides of Equality (\ref{hyperHproof}), since $\xi^{k-2}(x)\oplus\xi^{n-k+2}(w)=E$ by 2-hyperconvexity, we get
\begin{eqnarray}\label{hyperHproof2}
(\xi^{n-k+1}(y)\cap \xi^{k}(x))+\xi^{k-2}(x) +(\xi^{n-k+1}(t)\cap \xi^{k}(x))&=&\xi^{k}(x).
\end{eqnarray}
Adding $\xi^{n-k}(y)$ to both sides of  Equality (\ref{hyperHproof2}), we get
\begin{eqnarray}\label{hyperHproof3}
\xi^{n-k+1}(y)+\xi^{k-2}(x) +(\xi^{n-k+1}(t)\cap \xi^{k}(x))&=&E.
\end{eqnarray}
Now for dimension reasons, the above is direct and we get the following equality which is nothing else than Property (H), 
$$
\xi^{n-k+1}(y) \oplus\xi^{k-2}(x)
\oplus(\xi^{n-k+1}(t)\cap \xi^{k}(x)) =E.
$$
Conversely, let's assume Property (H), 
This implies after taking the intersection with $\xi^{k}(x)$ to
$$
(\xi^{n-k+1}(y)\cap \xi^{k}(x)) \oplus\xi^{k-2}(x)
\oplus(\xi^{n-k+1}(t)\cap \xi^{k}(x)) =\xi^{k}(x).
$$
Let denote $\pi$ the projection along $A=\xi^{k-2}(x)$ on $B=\xi^{n-k-2}(w)$.
Let $L(m)$ be the line $\xi^{n-k+1}(m)\cap \xi^{k}(x)$. The last assertion 
can be restated as
$$
L(y)\oplus L(t)\oplus A = \xi^{k}(x).
$$
Using $\pi$ we get
$$
\pi(L(t))\oplus\pi(L(y))=\pi(\xi^{k}(x))=\xi^{k}(x) \cap\xi^{n-k+2}(w),
$$
which is precisely Assertion (\ref{hyperHproof}), we wanted to show.
\qed

\subsubsection{Invariant subbundle}
Let's start with a definition just used in the next proof. A
subbundle $L$ of $E$ is said to be {\it invariant} if
\begin{displaymath}
\psi_t (L_x)= L_{\phi_t(x)}.
\end{displaymath}
Notice that no regularity is assumed on $L$.

\begin{lemma}\label{invar-sub}
  Let $L$ be an invariant subbundle of rank 1 of
  $E=W^1\oplus W^2$ equipped with a flag action. Assume
  that, for some auxiliary (continuous) metric $d$ on the bundle $\mathbb
  P ( E)$, we have
\begin{displaymath}
\exists\epsilon,\ \forall x \in M, \ d(L_x, W^1_x)>\epsilon>0.
\end{displaymath}
Then
\begin{displaymath}
\forall x \in M, \ L_x=W^2_x.
\end{displaymath}
\end{lemma}
We stress again that no regularity is assumed on $L$.

\proof The lemma is an immediate consequence of our
contraction Assumption (7).  Indeed if the lemma is not
true, for some $x$, $L_x\not= W^{2}_x$. Then, for some large
positive $s$ by our contraction Assumption (7)
$$
d(L_{\phi_{s}(x)},W^1_x)=d(\psi_{s}(L_x),W^{1}_x)<
\frac{\epsilon}{2}.
$$
\qed
\subsubsection{Proof of Lemma \ref{rank2}}
\proof By assumption, we have that for $y\not=z$,
$W^1_{x,y}\oplus W^1_{x,z}=E_x$.  Therefore the continuous
maps $J_x$ are injective. Since ${\mathbb P}(E_x)$ is of
dimension 1, the following limits exist (after a choice of
orientation on $\mathcal F$)
\begin{eqnarray*}
\lim_{y\rightarrow +\infty}J_{x}(y)&=&J^+_x\\
\lim_{y \rightarrow -\infty}J_{x}(y)&=&J^-_x
\end{eqnarray*}
Notice also that the bundles $J^{\pm}$ are flow invariant,
although not {\it a priori} continuous.  To finish proving
the lemma, we just have to show that
$$
J^+_x=J^-_x=W^{2}_x.
$$
We can assume that $\mathcal F$ is the orbit lamination
of a flow $\theta_{t}$.  Let's introduce the continuous
bundles $L^{\pm}_x=W^1_{x,\theta_{\pm 1}(x)}$. We notice
that for all $x$, $L^\pm_x\not=W^1_x$ by our assumption
(\ref{hyperconvex}). And furthermore, since the maps $J_x$
are monotone, $J^+_x$ and $W^1_x$ are not in the same
connected component of ${\mathcal P}(E_x)\setminus
\{L^-_x,L^+_x\}$ (and the same holds for $J^-_x$).  It
follows there exists $\epsilon >0$ such that
$$
d(J^\pm_x,W^1_x)\geq d(L^\pm_{x},W^1_x)\geq \epsilon.
$$
Lemma \ref{invar-sub} implies that for all $x$, $J^\pm_x=
W^{2}_x$. \qed

\subsection{Property (H) and hyperconvex bundles}
We recall that a $n$-Anosov representation, with limit curve
$\xi$, satisfies {\em Property (H)}, if for every triple of
distinct points $x$, $y$ and $z$ and integer $k$, we have
\begin{eqnarray*}
\xi^{k+1}(y)\oplus(\xi^{k+1}(z)\cap\xi^{n-k}(x))\oplus\xi^{n-k-2}(x)&=&E.
\end{eqnarray*}
We explain now various ways to check this property, and use 
its relation with hyperconvex bundles to  deduce an
important consequence in Proposition \ref{mainH}.
\subsubsection{Main property}\label{mainHdef}
Let $x$ and $z$ be two distinct points of
$\partial_{\infty}\grf$. Let
$$
G_{p,x,z}=\xi^{n-p+2}(z)\cap\xi^{p}(x).
$$
We consider the map $\mathcal Y_{p,x,z}$ defined by
\begin{displaymath}
\left\{
\begin{array}{rcl}
\partial_{\infty}\grf\setminus\{x\}&\rightarrow&\mathbb P (G_{p,x,z})\setminus \{\xi^{p-1}(x)\cap\xi^{n-p+2}(z)\}\\
y&\mapsto &(\xi^{n-p+1}(y)\oplus \xi^{p-2}(x))\cap G_{p,x,z}.
\end{array}
\right.
\end{displaymath}
Proposition \ref{preH} explains that this application is
well defined. Our main result in this paragraph is the
following Proposition
\begin{proposition}\label{mainH} Assume the representation satisfies Property (H). Then, the map $\mathcal Y_{p,x,z}$ is onto from $\partial_{\infty}\grf\setminus\{x\}$ to $\mathbb P (G_{p,x,z})\setminus \{\xi^{p-1}(x)\cap\xi^{n-p+2}(z)\}$.
\end{proposition}

\subsubsection{Construction of the map $\mathcal Y$}
We use the notation of the previous paragraph
\begin{proposition}\label{preH}
  We have
\begin{eqnarray}
\dim(\mathcal Y_{p,x,z}(y))&=&1,\label{preH2}\\
\mathcal Y_{p,x,z}(y)\oplus \big(\xi^{p-1}(x)\cap\xi^{n-p+2}(z)\big)&=&G_{p,x,z}.\label{preH3}
\end{eqnarray}
Furthermore, \label{equiH} Property (H) is equivalent to the
following assertion
\begin{eqnarray}
\mathcal Y_{p,x,z}(t)&\not=&\mathcal Y_{p,x,z}(y),\label{prehyp1}.
\end{eqnarray}
\end{proposition}
\proof From Identification (\ref{identif2}), we get that
$$
W^{1}_{(x,x_{0},w),(x,x_{0},y)}=\mathcal Y_{k,x,w}(y).
$$
All the above results follows from Paragraph
\ref{hbar}.\qed

\subsubsection{Proof of Proposition \ref{mainH}}
\proof From Proposition \ref{rkhc}, we know the bundles
$F_{k}=E^{k}/E^{k-2}$ are hyperconvex if the representation
satisfies Property (H). Furthermore, from Identification
(\ref{identif2}), we get that
$$
W^{1}_{(x,x_{0},w),(x,x_{0},y)}=\mathcal Y_{k,x,w}(y).
$$
Then Proposition \ref{mainH} follows from Lemma
\ref{rank2}.  \qed
\subsection{Proof of Lemma \ref{mainlemmacurve2}}
\subsubsection{Preliminary facts: case $l=1$.}
Let $\xi$ be the limit curve of the Anosov representation
$\rho$.  Let $(z,x,x_0,x_1,y)$ be distinct points of
$\partial_{\infty}\grf$.  Write $X=(z,x,x_{0},x_{1})$. We
introduce
\begin{eqnarray*}
U_{y,x_{0}}&=&(\xi^{n-k}(y)\oplus\xi^{k-1}(x_0)\cap\xi^{n-k+1}(z)),\\
Z_{x_{0},x_{1}}&=&(\xi^{k-1}(x_0)\oplus\xi^1(x_1))\cap\xi^{n-k+1}(z).
\end{eqnarray*}
We first notice that thanks to 2-hyperconvexity ({\em cf.}
Assertion (\ref{mlc0})) the sums in the definition of
$Z_{x_{0},x_{1}}$ and $U_{y,x_{0}}$ are indeed direct.  We
shall now prove
\begin{lemma}\label{preliml=1}
  If $\xi$ is a limit curve of a 3-hyperconvex
  quasi-Fuchsian representation, then
\begin{eqnarray*}
\dim(Z_{x_{0},x_{1}})&=&1\\
\dim(U_{y,x_{0}})&=&n-k.
\end{eqnarray*}
And furthermore,
\begin{eqnarray*}
Z_{x_{0},x_{1}}\oplus U_{y,x_{0}}&=&\xi^{n-k+1}(z).
\end{eqnarray*}
\end{lemma}
\proof Let
$$
C_y=\xi^{n-k}(y)\oplus\xi^{k-1}(x_0).
$$
Let $\pi$ be the projection on $A=\xi^{n-k+1}(z)$ along
$B=\xi^{k-1}(x_0)$.  Notice that $A\oplus B=E$ thanks to
2-hyperconvexity.  Recall that
$$
\pi(W)=(W+B)\cap A.
$$
In particular
\begin{eqnarray*}
\pi(C_y)&=&C_{y}\cap A=U_{y,x_{0}}\\
\pi( \xi^1(x_1))&=&Z_{x_{0},x_{1}}.
\end{eqnarray*}
We begin by computing the dimensions of $Z_{x_{0},x_{1}}$
and $U_{y,x_{0}}$. We first notice that
$$
\dim Z_{x_{0},x_{1}}=\dim(\pi(\xi^1(x_1))\leq 1.
$$
By 2-hyperconvexity, the next sum is direct
$$
\xi^1(x_1)+ \ \ \underbrace{\xi^{k-1}(x_0)}_B.
$$
Hence
$$
\xi^1(x_1)\not\subset B,
$$
and
$$
\dim(Z_{x_{0},x_{1}})=1.
$$
Let's consider now $U_{y,x_{0}}$. First, we know that
$$
\dim(C_y)=n-1.
$$
Finally, since the curve is assumed to be 3-hyperconvex,
$$
\underbrace{\xi^{n-k}(z)}_{\subset A}\oplus\ \ 
\underbrace{\xi^{k-1}(x_0)\oplus\xi^{1}(y)}_{\subset C_y}=E.
$$
Hence
$$
A+C_y=E,
$$
and
$$
A\not\subset C_y.
$$
Hence,
$$
\dim(U_{y,x_{0}})=\dim(\pi(C_y))=\dim(C_{y}\cap
A)=\dim(A)-1=n-k.
$$
Finally, since the curve is 3-hyperconvex,
\begin{eqnarray}
\underbrace{\big(\xi^{n-k}(y)\oplus\xi^{k-1}(x_{0})\big)}_{ C_{y}} \oplus \  \xi^1(x_1)&=&E.\label{mlc41}
\end{eqnarray}
Applying $\pi$ to both sides of Formula (\ref{mlc41}), we
finally get
$$
Z_{x_{0},x_{1}}+U_{y,x_{0}}=\xi^{n-k+1}(z).
$$
This concludes our proof.  \qed
\subsubsection{Proof of Lemma \ref{mainlemmacurve2}: case $l=1$}
We concentrate on the case $l=1$ of Lemma
\ref{mainlemmacurve2}. We prove
\begin{lemma}\label{mainlemmacurvel=1}
 Assume the limit curve $\xi$ of the Anosov representation
 $\rho$ is
\begin{itemize}
\item 3-hyperconvex
\item satisfies Property (H).
\end{itemize}
Then it is $(k,1)$-convergent for all $k$.
\end{lemma}
\proof We shall prove this Lemma by induction on $n-k$ and
use the following induction hypothesis.
\begin{eqnarray}
\lim_{(x_0,x_1)\rightarrow x}\big(\xi^k(x_0)\oplus\xi^1(x_1)\big)&=&\xi^{k+1}(x),\label{mlbc10}
\end{eqnarray}
We notice this is true for $k=n-1$ thanks to 2-hyperconvexity.  We want to prove
\begin{eqnarray}\label{mlcf1}
\lim_{(x_0,x_1)\rightarrow x}\big(\xi^{k-1}(x_0)\oplus\xi^1(x_1)\big)&=&\xi^{k}(x).
\end{eqnarray}
Recall that by 2-hyperconvexity, for $z\not=x$,
\begin{eqnarray}
\xi^{k-1}(x)\oplus\xi^{n-k+1}(z)=E.\label{mlbc11}
\end{eqnarray}
Let's introduce as in the previous paragraph
\begin{eqnarray*}
Z_{x_{0},x_{1}}&=&(\xi^{k-1}(x_0)\oplus\xi^1(x_1))\cap\xi^{n-k+1}(z),\\
U_{y,x_{0}}&=&(\xi^{n-k}(y)\oplus\xi^{k-1}(x_0))\cap\xi^{n-k+1}(z),\\
B&=&\xi^{n-k+1}(z),\\
G_{x,z}&=&B\cap\xi^{k+1}(x).
\end{eqnarray*}
Using Assertion (\ref{mlbc11}) and the notations of the
previous paragraph, since $\xi^{k-1}(x_0)$ converges to $\xi^{k-1}(x)$, we first notice Assertion (\ref{mlcf1}) would 
follow from
$$
\lim_{(x_0,x_1)\rightarrow
  x}Z_{x_{0},x_{1}}=\xi^{n-k+1}(z)\cap\xi^{k}(x).
$$
Our aim now is to prove this last assertion. We shall do
that by ``trapping" $Z_{x_{0},x_{1}}$ using $U_{y,x_{0}}$.

Let $V$ be a connected neighbourhood of $x$ homeomorphic to
the interval. Let's choose an orientation on $V$.  We shall
say $(x_0,x_1)$ tends to $x_+$, (resp. $x_-$) if $x_0>x_1$
(resp. $x_0<x_1$). Let $\alpha\in\{+,-\}$. Let
$\Lambda^\alpha_x$ be the set of value of accumulation of
$Z_{x_0,x_1}$ when $(x_0,x_1)$ tends to $x_\alpha$.  We note
that $\Lambda^\alpha_x$ is a connected subset a of ${\mathbb
  P}(B)$.  Notice that,
$$
Z_{x_{0},x_{1}}\subset( \xi^{k}(x_0)\oplus\xi^1(x_1))\cap
\xi^{n-k+1}(z).
$$
From Hypothesis (\ref{mlbc10}), it follows that the sets
$\Lambda^\alpha_x$ is actually a subset of the projective
space ${\mathbb P}(G_{x,z})$.  Therefore $\Lambda^\alpha_x$
is a closed interval of the 1-dimensional manifold ${\mathbb
  P}(G_{x,z})$.

We now choose an auxiliary metric $\langle,\rangle$ on $B$,
a unit vector $z_{x_{0},x_{1}}$ in $Z_{x_{0},x_{1}}$
continuous in $(x_{0},x_{1})$, a normal vector $u_{y,x_{0}}$
to $U_{y,x_{0}}$ continuous in $x_{0}$ and $y$. Then we have
(maybe after replacing $u$ by $-u$) for all $y$, for
$X=(x_0,x_1)$ close to $x$, thanks to Lemma \ref{preliml=1},
$$
\langle z_{x_{0},x_{1}},u_{y,x_{0}}\rangle>0.
$$
We consider now $\hat\Lambda^\alpha_x$ the set of value
of accumulation of $z_{x_{0},x_{1}}$ as $(x_{0},x_{1})$ goes
to $\alpha$.

Notice then, that for all $y$, for all $w$ in
$\hat\Lambda^\alpha_x$, we have
$$
\langle w,u_{y,x}\rangle\geq 0.
$$
Hence, for all $y, t$, $\hat\Lambda^\alpha_x$ is
contained in the closure of one connected component of
$$
G_{x,z}\setminus \big((U_{y,x}\cap G_{x,z}) \cup
(U_{t,x}\cap G_{x,z})\big).
$$
Since, by Definition \ref{mainHdef}
$$
U_{y,x}\cap G_{x,z} ={\mathcal Y}_{k+1,x,z}(y),
$$
It follows $\Lambda^\alpha_x$ is contained in the closure
of one connected component of
$$
{\mathbb P} (G_x)\setminus \{{\mathcal Y}_{k+1,x,z}(y),
{\mathcal Y}_{k+1,x,z}(t)\}.
$$
From Proposition \ref{mainH}, we know that the map
${\mathcal Y}_{k+1,x,z}$ is onto from
$\partial_{\infty}\grf\setminus\{x\}$ to $\mathbb P
(G_{k,x,z})\setminus \{\xi^{k}(x)\cap\xi^{n-k+1}(z)\}$.
Hence, we have that
$$
\Lambda^\alpha_x=\{\xi^{k}(x)\cap\xi^{n-k+1}(z) \}.
$$
Since this is true for all the ends $\alpha$, we obtain
the desired result. \qed

\subsubsection{Main Lemma: case $l>1$.}
The Main Lemma \ref{mainlemmacurve2} will follow by an
induction proved in Paragraph \ref{induction} from the next
statement combined with Lemma \ref{mainlemmacurvel=1}
\begin{lemma}\label{mainlemmacurve}
  Let $\xi$ be the limit curve of a quasi-Fuchsian
  representation. Let $k$ and $l$ be some integers such that
  $k+l\leq n$ and $l>2$.  We assume furthermore that the
  curve is
\begin{enumerate}
\item $(k,l-1)$-convergent,\label{mlc1}
\item $(k,l)$-convergent,\label{mlc2}
\item $(k-1,l-1)$-convergent, \label{mlc3}
\end{enumerate}
Then, the curve is $(k-1,l)$-convergent.
\end{lemma}

\subsubsection{Preliminary facts: case $l>1$}
We assume now that the limit curves $\xi$ satisfy the
hypothesis of Lemma \ref{mainlemmacurve}. That is we assume
the curve is
\begin{enumerate}
\item $(k,l-1)$-convergent,\label{mlc1b}
\item $(k,l)$-convergent,\label{mlc2b}
\item $(k-1,l-1)$-convergent, \label{mlc3b}
\end{enumerate}
Notice that since the curve is $(k,l)$-direct by Hypothesis  (\ref{mlc2b}) and Lemma \ref{drsl}, the next sum
is  direct for all $m$ with $m\leq k$,
\begin{equation}
\xi^{m}(x_0)+\xi^1(x_1)+\ldots+\xi^1(x_l)).\label{toto}
\end{equation}
Let $z$ be a point of   $\partial_\infty \grf$, $Y=(y_0,y_1,\ldots,y_{l-1}$ be a $l$-uple of cyclically ordered points of
 $\partial_\infty \grf\setminus\{z\}$. We denote by
\begin{eqnarray*}
C(z,Y)&=&\xi^{n-k-l}(z)\oplus\xi^{k}(y_0)\oplus\bigoplus_{i=1}^{i=l-1}\xi^1(x_i).
\end{eqnarray*}
The sum in the definition of $C(z,Y)$ is direct thanks to Hypothesis \ref{mainlemmacurve}.(\ref{mlc1}) and
Lemma \ref{drsl}. We also need to make some choices of orientation. Let
$I=\partial_{\infty}\grf\setminus\{z\}$. Let's choose an
orientation on $\xi^p(w)$ for all $p$ depending continuously
on $w$ in $I$. Let's also choose an arbitrary orientation on
$\xi^k(z)$ for all $k$. It follows that there exists a family of 1-forms $\alpha(z,Y)$ continuous in $Y$ such that
\begin{eqnarray*}
C(z,Y)&=&ker (\alpha(z,Y)).
\end{eqnarray*}

Let now $X=(z,x_0,x_1,\ldots,x_l)$ be a $l+2$-uple of
distinct points of $\partial_{\infty}\grf$ cyclically
oriented. Let 
\begin{eqnarray*}
X^+&=&(x_0,\ldots,x_{l-1})\\
X^-&=&(x_0,\ldots,x_{l-2},x_l).
\end{eqnarray*}
We introduce
\begin{eqnarray*}
U^+_{X}&=&C(z,X^+)\cap\xi^{n-k-l+2}(z),\\
U^-_{X}&=&C(z,X^-)\cap\xi^{n-k-l+2}(z),\\
Z_X&=&(\xi^{k-1}(x_0)\oplus\xi^1(x_1)\oplus\ldots\oplus\xi^1(x_l))\cap\xi^{n-k-l+2}(z).
\end{eqnarray*}
We first notice that thanks Assertion (\ref{toto}) the sum in the definition
of $Z_X$ is indeed direct. 

We shall now prove
\begin{proposition}
  If $\xi$ is a limit curve of a quasi-Fuchsian
  representation which satisfies the hypothesis of Lemma
  \ref{mainlemmacurve}, then
\begin{eqnarray*}
\dim(Z_X)&=&1\\
\dim(U^\pm_{X})&=&n-k-l+1.
\end{eqnarray*}
Furthermore,
\begin{eqnarray}
Z_X\oplus U^+_{X}&=&Z_X\oplus U^-_{X}=\xi^{n-k-l+2}(z),
\end{eqnarray}
and considering orientations 
\begin{eqnarray}
[\alpha(z,X^+)]\vert_{Z_X}=-[\alpha(z,X^-)]\vert_{Z_X}.\label{orienZU}
\end{eqnarray}
\end{proposition}
\proof Let
\begin{eqnarray*}
C^+&=&C(z,X^+)=\xi^{n-k-l}(z)\oplus\xi^{k}(x_0)\oplus\xi^1(x_1)\oplus\ldots\oplus\xi^1(x_{l-2})\oplus\xi^1(x_{l-1}),\\
B^+&=&\xi^{k-1}(x_0)\oplus\xi^1(x_1)\oplus\ldots\oplus\xi^1(x_{l-2})\oplus\xi^1(x_{l-1}),\\
B^-&=&\xi^{k-1}(x_0)\oplus\xi^1(x_1)\oplus\ldots\oplus\xi^1(x_{l-2})\oplus\xi^1(x_{l}),\\
C^-&=&C(z,X^-)=\xi^{n-k-l}(z)\oplus\xi^{k}(x_0)\oplus\xi^1(x_1)\oplus\ldots\oplus\xi^1(x_{l-2})\oplus\xi^1(x_{l}).
\end{eqnarray*}
Let $\pi^\pm$ be the projection on $A=\xi^{n-k-l+2}(z)$
along $B^\pm$.  Notice that $A\oplus B^\pm=E$ thanks to
Lemma \ref{drsl} and Hypothesis
\ref{mainlemmacurve}.(\ref{mlc3}). Recall that
$$
\pi^\pm(W)=(W+B^\pm)\cap A.
$$
In particular
\begin{eqnarray*}
\pi^\pm(C^\pm)=C^\pm\cap A&=&U^\pm_{X}\\
\pi^+(\xi^1(x_l))&=&Z_X,\\
\pi^-(\xi^1(x_{l-1}))&=&Z_X.
\end{eqnarray*}
We first compute the dimensions of $Z_X$ and $U^+_{X}$. We
first notice that
$$
\dim Z_X=\dim(\pi^+(\xi^1(x_l))\leq 1.
$$
Finally since the following sum is direct ({\it cf}
Hypothesis \ref{mainlemmacurve}.(\ref{mlc2}))
$$
\xi^k(x_0)+\xi^1(x_1)+\ldots+\xi^1(x_{l-1})
+ \ \xi^1(x_l),
$$
It follows the next one is direct
$$
\xi^1(x_l)+ \ \ 
\underbrace{\xi^{k-1}(x_0)+\xi^1(x_1)+\ldots+\xi^1(x_{l-1})}_{B^+}.
$$
Hence
$$
\xi^1(x_l)\not\subset B^+,
$$
and
$$
\dim(Z_X)=1.
$$
Let's consider now $U^+_{X}$, the proof for $U^-_X$ being
symmetric. First, we know that
$$
\dim(C^+)=n-1.
$$
According to Hypothesis \ref{mainlemmacurve}.(\ref{mlc3})
and Lemma \ref{drsl},
$$
\underbrace{\xi^{n-k-l+2}(z)}_A\oplus\ \ 
\underbrace{\xi^{k-1}(x_0)\oplus\xi^{1}(x_1)\oplus\ldots\oplus\xi^1(x_{l-1})}_{\subset
  C^+}=E.
$$
Hence
$$
A+C^+=E,
$$
and
$$
A\not\subset C^+.
$$
It follows that
$$
\dim(U^+_{X})=\dim(\pi^+(C^+))=\dim(C^+\cap
A)=\dim(A)-1=n-k-l+1.
$$
Finally, by Hypothesis \ref{mainlemmacurve}.(\ref{mlc2})
and Lemma \ref{drsl},
\begin{eqnarray}
\underbrace{\big(\xi^{n-k-l}(z)\oplus\xi^{k}(x_0)\oplus\xi^1(x_1)\oplus\ldots\oplus\xi^1(x_{l-1}\big)}_{C^+} \oplus \  \xi^1(x_l)&=&E.\label{mlc4}
\end{eqnarray}
Applying $\pi^+$ on both sides of Formula (\ref{mlc4}), we
finally get
$$
Z_X+U^+_{X}=\xi^{n-k-l+2}(z).
$$
The same holds for $U^-_X$. Finally, it remains to check
the orientations on $Z_X\oplus U^+_{X}$ and $Z_X\oplus
U^-_{X}$ are opposite. We shall denote by $\overline V$ the opposite of the oriented vector space $V$.

Since $(z,x_{0},\ldots,x_{l})$ are distinct and cyclically
oriented, there exists an arc $t\mapsto w_t$ joining $x_l$
to $x_{l-1}$, such that
$$
\forall t, w_t\not\in \{z,x_0,\ldots,x_{l-2}\}.
$$
Let
$$
B_t=\xi^{k-1}(x_0)\oplus\xi^1(x_1)\oplus\ldots\oplus\xi^1(x_{l-2})\oplus\xi^1(w_t).
$$
Notice that $B_t$, as before, satisfies
\begin{eqnarray}
B_t\oplus \xi^{n-k-l+2}(z)=E.\label{orien}
\end{eqnarray}
Let's choose an orientation on $E$ such that with respect to
the orientation
$$
E=C^+\oplus
\xi^1(x_l)=\ldots\oplus\xi^1(x_{l-1})\oplus\xi^1(x_l).
$$
Recall $B_{t}$ is oriented. We choose now the orientation
on $\xi^{n-k-l+2}(z)$ compatible with Equation
(\ref{orien}). It follows that considered as oriented space
we have
$$
U^+_X\oplus Z_X=\xi^{n-k-l+2}(z).
$$
Conversely, since
$$
E=\overline{C^-\oplus
  \xi^1(x_{l-1}))}=\overline{\ldots\oplus\xi^1(x_{l})\oplus\xi^1(x_{l-1})}.
$$
We obtain that
$$
\overline{U^+_X\oplus Z_X}=\xi^{n-k-l+2}(z).
$$
\qed
\subsubsection{Proof of Lemma \ref{mainlemmacurve}}
Let's first state the following elementary lemma.
\begin{lemma}\label{elemorien}
  Let $E$ be a vector space, Let $\{L_n\}_{n\in\mathbb N}$
  be a sequence of oriented lines converging to an oriented
  line $L_\infty$.  Let $\{P_n\}_{n\in\mathbb N}$ be a
  sequence of oriented hyperplanes converging to an oriented
  hyperplane $P_\infty$.  Assume that the following sums are
  direct and with opposite orientations
\begin{eqnarray*}
L_n\oplus P_{2n}=\overline{L_n\oplus P_{2n+1}}.
\end{eqnarray*}
Then $L_\infty\subset P_\infty$.
\end{lemma}
\proof Let's choose an auxiliary metric $g$ on $E$. Let
$u_n$ be the positive unit vector in $L_n$ and $v_n$ be the
normal unit vector to $P_n$. From the hypothesis, we get
that
$$
g(u_n,v_{2n}).g(u_n,v_{2n+1})<0.
$$
Therefore, by passing to the limit we obtain that
$g(u_\infty,v_\infty)=0$. \qed \vskip 1truecm

Let's now proceed to the main proof.  We shall always assume
that
$$
(z,x_0,x_1,\ldots,x_l)
$$
are distinct and cyclically positively oriented. We
choose as before an orientation on $\xi^p(w)$ depending
continuously on $w$ in
$\partial_{\infty}\grf\setminus\{z\}$. Here are the
hypothesis we shall assume.
\begin{eqnarray}
\lim_{(x_0,\ldots,x_l)\rightarrow x}\big(\xi^k(x_0)\oplus\xi^1(x_1)\oplus\ldots\oplus\xi^1(x_l)\big)&=&\xi^{k+l}(x),\label{mlbc1}\\
\lim_{(x_0,\ldots,x_l)\rightarrow x}\big(\xi^k(x_0)\oplus\xi^1(x_1)\oplus\ldots\oplus\xi^1(x_{l-1})\big)&=&\xi^{k+l-1}(x),\label{mlbc2}\\
\lim_{(x_0,\ldots,x_{l-2})\rightarrow x}\big(\xi^{k-1}(x_0)\oplus\xi^1(x_1)\oplus\ldots\oplus\xi^1(x_{l-1})\big)&=&\xi^{k+l-2}(x).\label{mlbc3}
\end{eqnarray}
We can actually assume the limit in Assertion (\ref{mlbc2})
is a limit as oriented vector spaces.  We want to prove
\begin{eqnarray}\label{mlcf}
\lim_{(x_0,\ldots,x_l)\rightarrow x}\big(\xi^{k-1}(x_0)\oplus\xi^1(x_1)\oplus\ldots\oplus\xi^1(x_l)\big)&=&\xi^{k+l-1}(x).
\end{eqnarray}
We shall use the notations and results of the preceding
paragraph.

It follows that using Hypothesis (\ref{mlbc3}), Assertion
(\ref{mlcf}) reduces to
$$
\lim_{(x_0,\ldots,x_l)\rightarrow
  x}Z_{X}=\xi^{n-k-l+2}(z)\cap\xi^{k+l-1}(x).
$$
Our aim now is to prove this last assertion. We shall do
that by ``trapping" $Z_X$ using $U^\pm_{X}$.

Let $\Lambda_x$ be the set of value of accumulation of
$Z_{X}$ when $(x_0,\ldots,x_{l})$ tends to $x$.  We note
that $\Lambda_x$ is a subset of ${\mathbb
  P}(\xi^{n-k-l+2}(z))$.  Recall that
$$
Z_{X}\subset
\xi^{k}(x_0)\oplus\xi^1(x_1)\oplus\ldots\oplus\xi^1(x_l).
$$
From Hypothesis (\ref{mlbc1}), we finally get the set
$\Lambda_x$ is actually a subset of the projective space
${\mathbb P}(W)$, for $W=\xi^{n-k-l+2}(z)\cap\xi^{k+l}(x)$.

Finally, from Hypothesis (\ref{mlbc1})
$$
\lim_{X\rightarrow x}(U^\pm_X)=\underbrace
{(\xi^{n-k-l}(z)\oplus\xi^{k+l-1}(x))\cap\xi^{n-k-l+2}(z)}_D.
$$
Notice that $\xi^{n-k-l+2}(z)+ \xi^{k+1}(x)=E$ by 2-hyperconvexity hence $D$ is indeed a hyperplane of  $\xi^{n-k-l+2}(z)$.
We can choose the orientation on $D$ such the limit is to
be considered for oriented vector spaces. From Equation
(\ref{orienZU}) and Lemma \ref{elemorien}, we obtain that
$$
\Lambda_x\subset\mathbb P (D).
$$
The conclusion of the proof follows from
\begin{eqnarray*}
D\cap W&=&\big(\xi^{n-k-l}(z)\oplus\xi^{k+l-1}(x)\big)\cap\xi^{n-k-l+2}(z)\cap\xi^{k+l}(x)\\
&=&\xi^{k+l-1}(x)\cap\xi^{n-k-l+2}(z).
\end{eqnarray*}
\qed
\subsubsection{Final induction}\label{induction}
\proof It remains to prove Main Lemma \ref{mainlemmacurve2},
using Lemma \ref{mainlemmacurvel=1} and Lemma
\ref{mainlemmacurve}. This is done by induction. We say a
limit curve is {\em $l$-superconvergent}, if it is
$(k,l)$-convergent for all $k$.

From Lemma \ref{mainlemmacurvel=1}, the curve if
$1$-superconvergent. We assume by induction the curve is
$l-1$-superconvergent.  From Lemma \ref{mainlemmacurve} and
an easy induction, to prove that the curve is
$l$-superconvergent, it suffices to show that it is
$(n-l,l)$-convergent. But to be $(n-l,l)$-convergent just
means that the following sum is direct
$$
\xi^{n-l}(x_{0})+\xi^{1}(x_{1})+\ldots+\xi^{1}(x_{l})=E.
$$
But the fact this sum is direct follows from the fact the
curve is $(1,l-1)$-convergent, hence $(1,l-1)$-direct by
Lemma \ref{drsl}. All these conditions are guaranteed by the
induction assumption. \qed

\section{Anosov representations, Pro\-per\-ty (H) and  3-hyperconvexity }\label{5}
We now clarify some relations between Property (H),
3-hyperconvexity, and Anosov representations.  We first say
a representation is {\em $S$-irreducible} if its restriction
to all finite index subgroups is irreducible.  By Lemma
\ref{corohitchin}, every representation in Hitchin's
component is $S$-irreducible.  We shall denote
\begin{itemize}
\item $\mathcal A$ (resp. $\qf$) the space of $n$-Anosov
  $S$-irreducible representations (resp. quasi-Fuchsian
  representation),
\item $\mathcal A_H$ (resp. $\qf_H$) the space of
  $S$-irreducible Anosov (resp. quasi-Fuchsian)
  representations satisfying Property (H)
\item $\mathcal A_{3}$ (resp. $\qf_{3}$) be the set of
  $S$-irreducible Anosov (resp. quasi-Fuchsian) which are
  3-hyperconvex.
\end{itemize} 
We summarise in the next Proposition the results of this
section.
\begin{proposition}\label{qfh}$\mathcal A_{3}$ is open in $\mathcal A$.
  $\mathcal A_H$ is a connected subset of $\mathcal A$.
  Furthermore $\qf_H=\qf$, and every Fuchsian representation
  is 3-hyperconvex.
\end{proposition}

The proof of Proposition \ref{qfh} follows the following
path: we will prove in the next paragraph that $\mathcal
A_H$ and $\mathcal A_{3}$ are open in $\mathcal A$, and in
the next one that $\mathcal A_H$ is closed in $\mathcal A$;
finally we prove that every Fuchsian representation is
3-hyperconvex and satisfies Property (H). This will complete
the proof of Proposition \ref{qfh}.
\subsection{Open}
We first notice.
\begin{proposition}
  The sets $\mathcal A_H$ and $\mathcal A_{3}$ are open in
  $\mathcal A$.
\end{proposition}
\proof This follows at once from the fact
$(\partial_{\infty}\grf^3\setminus\Delta)/\grf$ is compact
and that the conditions defining 3-hyperconvexity and
Property (H) are open in the corresponding product of flag
manifolds. \qed
\subsection{Closed}

The aim of this paragraph is to prove the following
assertion.

\begin{proposition}
  The set $\mathcal A_H$ is closed in $\mathcal A$.
\end{proposition}
\proof Let consider $\rho$ a $S$-irreducible Anosov
representation limit of representations in $\mathcal A_H$.
Let $\xi$ be the associated curve, and $\mathcal Y =\mathcal
Y_{k,x,z}$, the associated map defined in Paragraph
\ref{mainHdef}. By Proposition \ref{equiH}, we wish to prove
that $\mathcal Y$ is injective. Since $\mathcal Y$ is a
limit of continuous injective maps of a 1-dimensional
manifold into another, $\mathcal Y$ is monotone.

Therefore, if $\mathcal Y$ fails to be injective, there is
an open set $U$ in $\partial_{\infty}\grf$ on which it is
constant. We will prove this last assertion leads to a
contradiction.  Indeed we will prove the next assertion that
contradicts Lemma \ref{conseqirredu}.

\vskip 0.5truecm {\it Assertion.} There exists some
$(n-k-1)$-plane $A$, such that
\begin{equation}
\forall y\in U,\ \dim\xi^{k+1}(y)\cap A\geq 1\label{qfh1}.
\end{equation}
Notice first that
\begin{eqnarray*}
\mathcal{Y}(y)&=&(\xi^{k+1}(y)\oplus \xi^{n-k-2}(z))\cap G_{k,x,z}\\
&=&(\xi^{k+1}(y)\oplus \xi^{n-k-2}(z))\cap\xi^{n-k}(z)\cap \xi^{k+2}(x)\\
&=&((\xi^{k+1}(y)\cap\xi^{n-k}(z))\oplus \xi^{n-k-2}(z))\cap \xi^{k+2}(x).
\end{eqnarray*}
Since $\xi^{k+2}(x)$ is a supplementary of $\xi^{n-k-2}(z)$,
$$
\mathcal{Y}(y)=\mathcal{Y}(t),
$$
implies
$$
P(y)=P(t),
$$
where
$$
P(y)= (\xi^{k+1}(y)\cap\xi^{n-k}(z))\oplus
\xi^{n-k-2}(z).
$$
Notice that $P(y)$ has dimension $n-k-1$. As a
conclusion, we get Assertion (\ref{qfh1}).  \qed
\subsection{Back to Fuchsian representations}
We prove now that $\qf_H$ is not empty. More specifically
\begin{lemma}
  Every Fuchsian representation satisfies Property (H).
\end{lemma}

\proof First we notice that for every distinct points $x$
and $z$ the follwing sum is direct
$$
\xi^{k+1}(z)+\xi^{n-k-2}(x),
$$
hence, next one is also direct,
$$
(\xi^{k+1}(z)\cap\xi^{n-k}(x))+\xi^{n-k-2}(x) =P(z,x).
$$
It follows that if a Fuchsian representation does not
satisfy Property (H), then there exists a triple of distinct
points $(x,y,z)$ such that
\begin{eqnarray*}
\dim(\xi^{k+1}(y) \cap P(z,x))>0.
\end{eqnarray*}
In the case of a Fuchsian representation, the limit curve is
the Veronese embedding and is equivariant under the whole
action of $SL(2,\mathbb R)$. Since $SL(2,\mathbb R)$ acts
transitively on the set of triple of distinct points, we
obtain there exists a $n-k-1$- plane $P$ (namely $P(z,x)$
for some $x$ and $z$) a such that for every $y$,
\begin{eqnarray*}
\dim(\xi^{k+1}(y) \cap P)>0.
\end{eqnarray*}
It follows that for every there exist a $k+1$-plane $Q$ such
that $A$ in $SL(2,\mathbb R)$,
\begin{eqnarray*}
\dim(\rho(A)Q\cap P)>0.
\end{eqnarray*}
This last assertion contradicts Proposition
\ref{algebirredu}.  \qed \vskip 1truecm Actually, one could
prove the previous proposition by an explicit computation.
Indeed, if we identify $\partial_\infty\grf\setminus\{x\}$
with $\mathbb R\mathbb P^1\setminus\{\infty\}=\mathbb R$,
and use the fact the irreducible representation of
$SL(2,\mathbb R)$ of dimension $n$ is the representation on
homogeneous polynomials of degree $n-1$ in variables $t$ and
$s$, one sees that
$$
\xi^k (x)=\{P(s,t)/ \exists Q {\hbox{ such that
  }}=(s+tx)^{n-k}Q(s,t)\}.
$$
Then it is an exercise (left to the reader) to prove the
previous proposition along these lines.

Similar arguments show
\begin{proposition}
  Every Fuchsian representation is 3-hyperconvex.
\end{proposition}

\section{Closedness}\label{8}

Our aim in this section is to prove the following result

\begin{lemma}\label{mainclosed}
  The set
  $$
  \tilde{\mathcal A}=\mathcal A_{3}\cap\mathcal A_{H}
  $$
  of 3-hyperconvex Anosov representations satisfying
  Property (H) is closed in the space of $S$-irreducible
  representations.
\end{lemma}
This Lemma will be deduced from Lemma \ref{genhypercurve}.
We first show that as corollaries, we obtain our Theorem \ref{maincurve}.
\subsection{Proof of Theorems \ref{maincurve}}
\subsubsection{Theorem  \ref{maincurve}}\label{proof1}

\proof We just have to put the previous statements in the
correct order. We first know by Lemma \ref{corohitchin} that
every representation in Hitchin's component is
$S$-irreducible.  By Proposition \ref{qfh},
$\qf_{H}\cap\qf_{3}$ is open in $\qf$, hence in Hitchin's
component by Lemma \ref{Thurslok}. It is non empty by
Proposition \ref{qfh} again, since it contains all Fuchsian
representations.  It is closed by Lemma \ref{mainclosed}.
Hence $\qf_{H}\cap\qf_{3}$ is equal to the all of Hitchin's
component.  Let $\rho$ be a representation in this
component. Let
$$
\xi=(\xi^1,\xi^2,\ldots,\xi^{n-1})
$$
be its limit curve.  By Corollary \ref{mainlemmacurve3},
we know that if
$\rho\in\qf_{H}\cap\qf_{3}=Rep_{H}(\grf,PSL(n,\mathbb R))$,
$\xi^{1}$ is an hyperconvex Frenet curve and $\xi$ is its
osculating flag \qed
%


\subsection{Convergence of limit curves}

Our aim is to prove the following lemma,

\begin{lemma}\label{genhypercurve}
  Let $\{\rho_m\}_{m\in \mathbb N}$ be a sequence of Anosov
  representations satisfying Property (H) converging to a
  $S$-irreducible representation $\rho$. Let
  $\xi_m=(\xi_m^1,\ldots,\xi^{n-1}_m)$ be the limit curve of
  $\rho_m$.  Then, there exists
\begin{itemize} 
\item a sequence of homeomorphisms $\phi_{m}$ of $S^{1}$
  with $\partial_{\infty}\grf$,
\item a monotone map $\pi$ from $S^{1}$ to
  itself,
\item an injective map $\hat\xi^1$
  from $S^1$ to $\mathbb P (E)$,
\item an injective left continuous orientation preserving map $\phi_0$ from $\partial_\infty\grf$ to $S^1$
\end{itemize}
such that 
\begin{itemize}
\item after extracting a subsequence the mappings
$\{\xi_m^1\circ\phi_{m}\}_{m\in \mathbb N}$ converges to
$\hat\xi^1\circ\pi$,
\item the map $\hat\xi^1\circ\phi_0$ is $\rho$-equivariant and $*$-hyperconvex.
\end{itemize}
\end{lemma}

We explain first this Lemma implies Lemma \ref{mainclosed} :
indeed, thanks to Theorem \ref{preservehyper} the limit
representation is Anosov and it satisfies Property (H)
thanks to Proposition \ref{qfh}.  The proof of the Lemma by
itself follows several steps which we describe now shortly
using the notations and the hypothesis of the Lemma.
\begin{enumerate}
\item {\em Convergence of the images (Proposition
    \ref{coi})}: there exists a sequence of homeomorphisms
  $\phi_m$ of $S^{1}$ with $\partial_{\infty}\grf$ such that
  $\{\xi_m^1\circ \phi_m\}_{m\in \mathbb N}$ converges to a
  rectifiable curve $\xi^1$
\item {\em Preliminary facts}: we prove lemmas of
  independent interest concerning rectifiable curves
  invariant under actions of groups.
\item {\em The limit and the boundary at infinity}: this is
  the core of the proof, in particular Lemma
  \ref{genhypercurve} is a consequence of Proposition
  \ref{lbi}. We basically prove that
  $\xi^{1}=\hat\xi^{1}\circ\pi$ where $\hat\xi^{1}$ is
 $*$-hyperconvex, $\rho$-equivariant, and $\pi$ is monotone
  from $S^{1}$ to $\partial_{\infty}\grf$.
\end{enumerate} 
From now on, we use the notation of the Lemma. That is we
consider
\begin{itemize}
\item a sequence $\{\rho_m\}_{m\in \mathbb N}$ of
  3-hyperconvex Anosov representations satisfying Property
  (H) converging to a $S$-irreducible representation $\rho$,
\item $\xi_m=(\xi_m^1,\ldots,\xi^{n-1}_m)$, the limit
  curve of $\rho_m$.
\end{itemize}

\subsection{Convergence of the images}
\begin{proposition}\label{coi}
  After passing to a subsequence, there exists a sequence of
  homeomorphisms $\phi_m$ of $S^{1}$ with
  $\partial_{\infty}\grf$such that $\{\xi_m^1\circ
  \phi_m\}_{m\in \mathbb N}$ converges to a rectifiable
  curve $\xi^1$.
\end{proposition}

\proof By Ascoli Theorem, it suffices to show that there
exists a constant $B$ such that $\xi^1_m$ is rectifiable and
with length bounded by $B$. But this follows from the next
remark: if $c$ is a curve in $\mathbb P(E)$, then
$$
{\hbox{length}}(c)\leq \int_{P(E^*)}\sharp (c\cap
P)d\mu(P).
$$
In our case, by Lemma \ref{mainlemmacurve2},
$\xi^{1}_{m}$ is hyperconvex, hence we get
$$
{\hbox{length}}(\xi^{1}_{m})\leq (dim(E)-1)\mu(P(E^*)).
$$
\qed
\subsection{Preliminary facts}

\subsubsection{Wormlike}\label{lamb}
Let $Z$ be a subset of $\mathbb P (E)$. We define $\langle
Z\rangle$ to be the vector subspace generated by all the
elements of $Z$:
$$
\langle Z\rangle=\sum_{u\in Z}u.
$$
Finally, assume $\Gamma$ acts on $S^1$. Let $\rho$ be a faithful representation of $\Gamma$ in $SL(E)$. Let $\xi$ be a $\rho$-equivariant injective map from  $S^1$ to $\mathbb P (E)$. Let
 $$
 \Gamma_{\mathbb R}=\{\gamma\in\Gamma, \gamma\not=id, \rho(\gamma) \hbox{ is diagonalisable over } \mathbb R \}.
 $$
 For every $\gamma$ in $\Gamma_{\mathbb R}$, let $Fix(\gamma)$ (resp. $Fix^+(\gamma)$) be the set of (resp. attractive) fixed points of $\rho(\gamma)$ in $\mathbb P (E)$. We define
$$
\Lambda_{\xi,\rho,\Gamma}=\{a\in S^1/ \exists \gamma\in\Gamma_{\mathbb R}, \xi(a) \in Fix^+(\gamma)\}.
$$  

We prove the following lemma.

\begin{lemma}\label{localworm}. 
  Let $\Gamma$ be a 
  group acting on $S^1$  by orientation preserving homeomorphisms. Let $\rho$ be a $S$-irreducible
  representation of $\grf$ in $SL(E)$.   Assume $\xi$ is a
  $\rho$-equivariant rectifiable  injective map from $S^1$ to $\mathbb P(E)$ with finite length. Then
  \begin{itemize}
  \item[$(i)$] $\xi(S^1)$ is not included in a finite union of proper vector subspaces of $E$,
  \item[$(ii)$] For every  $\gamma$ in $\Gamma_{\mathbb R}$, $S^1\setminus Fix(\gamma)$, has finitely many connected components.
  \item[$(iii)$] For every $\gamma$ in $\Gamma_{\mathbb R}$, there exists a unique  $\gamma^+$ in $S^1$ such that  $\xi(\gamma^+)\in Fix^+(\gamma)$. In particular, $\Lambda_{\xi,\rho,\Gamma}$ is not empty if $\Gamma_{\mathbb R}$ is not empty.
  \item[$(iv)$] Let $U$ be a neighbourhood of a point $c^+$ in $\Lambda_{\xi,\rho,\Gamma}$. Then 
  $$
  \xi(U)\not\subset P_0\cup P_1,
  $$
  for any proper vector subspaces $P_0$ and $P_1$ of $E$.
  In particular,
 $\langle \xi (U)\rangle=E$.
 \end{itemize}
\end{lemma} 
\proof  Statement $(i)$ is a consequence of the fact that the connected component of the identity of the Zariski closure of $\rho(\Gamma)$ is irreducible. Indeed, let $E_1,\ldots,E_p$ be proper vector subspaces such that
$$
\xi(S^1)\subset E_1\cup\ldots\cup E_p.
$$
We can as well assume that $<\xi(S^1)\cap E_i>=E_i$. It follows that for every $\gamma$ in $\Gamma$, one has 
$$
\gamma(E_i)\subset E_1\cup\ldots\cup E_p.
$$
The same property holds for $\gamma$ in the Zariski closure $H$ of $\rho(\gamma)$. Let $E_k$, such that $\dim(E_k)=\sup_i (\dim(E_i))$. Then, for every element $g$ in $H$ close to the identity, $g(E_k)=E_k$. It follows the identity component of $H$ preserves $E_k$, hence is not irreducible. This is the contradiction.

Let's now describe  the action  on $\mathbb P(E)$ of an element $f$ of $SL(E)$ diagonalisable over $\mathbb R$. These are elementary facts whose proofs are left to the reader.
\begin{itemize}
\item[(a)] The stable
  manifold $W$ of a fixed point $z$ of $f$ in $\mathbb P(E)$
  is described in the following way.  There exists a vector
  subspace $\tilde W$ of $E$, such that $W$ is an open set
  in $\mathbb P(\tilde W)$. Furthermore $W$ is open in $\mathbb P(E)$, if and only if $z$ is an attractive fixed point.
\item[(b)] Every closed invariant set of $f$ contains a fixed point.
\item[(c)] If $x$ is such that $f^n(x)$ converges to $a$ and $f^{-n}(x)$ converges to $b$ when $n$ goes to infinity with $a\not=b$, then $a$ and $b$ belong to different connected components of the space of fixed points of $f$.
\item[(d)] $f$ has at most one attractive fixed point.
\end{itemize}

We can now prove $(ii)$. Let $I=]\alpha,\beta[$ be a connected component of $S^1\setminus Fix(\gamma)$. Then $I$ is fixed by $\gamma$ since $\gamma$ is orientation preserving. Furthermore, by $(c)$, $\xi(\alpha)$ and $\xi(\beta)$  belong to different connected components of the space of fixed points of $\rho(\gamma)$. Let ${\mathcal W}$ be the set of connected components of $Fix(\rho(\gamma))$.  It follows that
$$
\hbox{length}(\xi(I))\geq \epsilon_0=\inf_{A,B\in\mathcal W, A\not=B} d(A,B).
$$
Since $\xi(S^1)$ has finite length by hypothesis, we deduce $(ii)$.

Let's proceed to $(iii)$. Assume that $Fix^+(\gamma)$ is empty. Notice that every connected component $]\alpha,\beta[$ of $S^1\setminus Fix(\gamma)$ is mapped to the stable manifold of $\xi(\alpha)$ and unstable manifold of $\xi(\beta)$ (after a choice of orientation). Therefore by $(a)$, if $Fix^+(\gamma)$ is empty, then $\xi(I)$ lies in a proper subspace of $E$. Since $\rho(\gamma)\not=id$, $\xi(Fix(\gamma))$ lies in a finite union of proper vector subspace. It follows that $\xi(S^1)$ lies in a finite (by $(ii)$) union of proper vector subspaces of $E$ and this contradicts $(i)$. Uniqueness follows from $(d)$ and the injectivity of $\xi$.

We shall now prove $(iv)$ by similar arguments. Let $c^+$ be the point which is mapped to the attractive fixed point of $\rho(\gamma)$, with $\gamma\in\Gamma_{\mathbb R}$. Write the finite decomposition in connected components
$$
S^1\setminus Fix(\gamma)= \bigsqcup_{i}V_i.
$$
By convention, we assume that $V_1$ and $V_2$ have $c^+$ in their closure. Let $i\geq 3$. The the closure of $V_i$ contains an element  which is mapped by $\xi$ to a fixed point $c$ which is neither attractive nor repulsive. The sets $V_i$ lie in the stable (or unstable) manifold of $c$. The same holds for $\xi (V_i)$. It follows from $(a)$ that $\xi (V_i)$  lies in a proper subspace $E_i$ of $E$, for $i\geq 3$.

Assume now  that there exist  a neighbourhood $U$ of $c^+$, two proper vector subspaces $P_0$ and $P_1$ of $E$, 
such that 
$$
\xi (U)\subset P_0\cup P_1.
$$
We choose $U$ small enough so that $U\subset \gamma^{-1}(U)$. Then, if $Q_i$  are limits of $\rho(\gamma^{n_q})(P_i)$ for some subsequence $n_q$, we get 
$$
\bigcup_{n\in \mathbb N} \xi (\gamma^{-n}(U))\subset Q_0\cup Q_1.
$$
But 
$$
\bigcup_{n\in \mathbb N}\gamma^{-n}(U)=\overline{V_1\cup V_2}.
$$
It follows that $\xi(S^1)$ lies in the union of $E_i$ for $i\geq3$, $Fix(\rho(\gamma))$ and $Q_0\cup Q_1$, hence the contradiction by $(i)$.
\qed

\subsubsection{Weak worm}
We prove now a weak version of the previous Lemma
\begin{lemma}\label{weakworm}
  Let $\xi$ be a rectifiable map parametrised by arc length
  from $S^1$ to $\mathbb P(E)$.  Let $\rho$ be an
  representation of $\grf$ in $SL(E)$. Assume $\rho$ is
  $S$-irreducible. Assume $\xi(S^{1})$ is
  $\rho(\grf)$-invariant.  Let $x,y$ be two distinct points
  of $S^{1}$, then one of the connected component $I$ of
  $S^{1}\setminus \{x,y\}$ satisfies
  $\langle\xi(I)\rangle=E$.
\end{lemma} 
The main point here is that $\xi$ is not assumed to be
injective. If it were, it would be an homeomorphism, and we
would have deduced an action of $\grf$ on $S^{1}$ such that
$\xi$ is $\rho$-equivariant and we could apply Lemma
\ref{localworm}.
 
\proof If both connected components $I_{0}$ and $I_{1}$ of
$S^{1}\setminus \{x,y\}$ satisfy
$$
\langle\xi(I_{i})\rangle=P_{i}\subsetneq E.
$$
Then $\xi(S^1)\subset P_0\cup P_1$,
 hence, $\rho$ would not be $S$-irreducible, by the same argument used in the proof of $(i)$ of the previous lemma, which also apply in this more general context.\qed

\subsection{The limit and the boundary at infinity}

From Proposition \ref{coi}, we can as well assume that
$\{\xi^1_m\}_{m\in\mathbb N}$ with the arc-length
parametrisation converges. Let in particular $\xi^1$ be the
limit. {\it A priori}, by using the arc-length
parametrisation, we have lost control over the action of
$\grf$. We just know that $\xi^1(S^{1})$ is globally
invariant by $\rho(\grf)$. Our aim now is to show this
action is semi-conjugate to the action of $\grf$ on
$\partial_{\infty}\grf$.

We begin by replacing $\xi^{1}$ by its arc-length parametrisation
$\hat\xi^{1}$ so that we have $\xi^{1}=\hat\xi^{1}\circ\pi$
with $\pi$ monotone.

We wish to prove

\begin{proposition}\label{lbi}
  There is an
  injective left continuous map preserving the orientation $\varphi_{0}$ from $\partial_{\infty}\grf$ to
  $S^{1}$, such that $\hat\xi^1\circ\varphi_{0}$ is $\rho$
  equivariant and $*$-hyperconvex.
\end{proposition}
Notice that Proposition \ref{lbi} implies Lemma
\ref{genhypercurve}.  The proof falls is several steps, we
prove
\begin{enumerate}
\item $\hat\xi^1$ is ``hyperconvex" when restricted to  a certain (non empty) subset $\Lambda$: Proposition \ref{bil};
\item finally, we prove Lemma \ref{northsouth} which,
  combined with the propositions of the previous section,
  implies Proposition \ref{lbi}.
\end{enumerate}

\subsubsection{$\Lambda$-Hyperconvexity}

We are going to prove the following two related propositions
\begin{proposition}\label{bil} 
The map $\hat\xi^1$ is injective.
\end{proposition}
As a consequence, it is an homeomorphism onto its image, and we deduce there exists an action of $\grf$ by homeomorphisms on $S^{1}$ such that $\hat\xi^1$ is $\rho$-equivariant. Let $\Gamma_0$ the normal subgroup of index 2 of orientation preserving elements of $\grf$. 
Let 
$$
\Lambda=\Lambda_{\xi,\rho,\Gamma_0}.
$$
Notice that
$\Lambda$ is $\grf$ invariant. We shall also prove.
\begin{proposition}\label{bilhyp}
For any $n$-uple of distinct points 
$(x_1,\ldots,x_n)$ of distinct points of the closed set $\overline\Lambda$ the following sum is direct
$$
\sum_{i=1}^{i=n}\hat\xi^1(x_i).
$$
\end{proposition} 
\proof Let's write
$\tilde\xi^{1}_{m}=\xi^{1}_{m}\circ\phi_{m}$, so that $$
\lim_{m\rightarrow\infty}\tilde\xi^{1}_{m}=\xi^{1}=\hat\xi^{1}\circ\pi.
$$
Let's prove the propositions. We split the proof in two parts, which are going to use very similar ideas. 

{\em Injectivity : proof of Proposition \ref{bil}.} First we want to prove the map is injective. Assume therefore that $\hat\xi^1(y)=\hat\xi^1(z)$. Since  $\rho$ is S-irreducible, by Lemma \ref{weakworm}, one of the connected component $J$ of
  $\partial_{\infty}\grf\setminus\{y,z\}$ is such that
  $\langle \xi(J)\rangle=E$.

We can therefore  find $n$
points $(x_{1},\ldots,x_{n})$ in $J$ such that the
following sums are direct
\begin{eqnarray}
\hat\xi^1(x_1)+\ldots+\hat\xi^1(x_n)&=&E,\cr
\forall i,\ \ \  \sum_{i\not=j}\hat\xi^1(x_j)+\hat\xi^1(y)&=&
E\label{goodsumdir}
\end{eqnarray}
Let $I$ be an interval containing $y$ and $z$ and none of the $x_i$.
For any of the points $t\in\{y,z,x_1,\ldots,x_n\}$, we denote by $\dot t$ a
point such that $\pi (\dot t)=t$.  For any distinct integers
$i,j$, let's write
$$
W_{ij}=(\dot x_i,\dot x_j), \ \ Y_{ij}=(\ldots,\dot
x_l,\ldots)_{l\not\in\{i,j\}}.
$$
We can as in Section \ref{hyperincrease0} consider the
maps $F^m_{ij}$ defined by
\begin{eqnarray*}
\mapping{ I}
{\mathbb P \big(\tilde\xi_{m}^{(2)}(W_{ij})\big)\setminus
\{\tilde\xi^1_m(\dot x_i)\}}
{t}
{\big(\tilde\xi_{m}^{(n-2)}(Y_{ij})\oplus\tilde\xi^1_{m}(t)
\big)\cap\tilde\xi_{m}^{(2)}(W_{ij})}.
\end{eqnarray*}
By Assertion (\ref{goodsumdir}), we
obtain that
$$
\lim_{m\rightarrow\infty} (F^m_{ij}(\dot y))
=\big(\bigoplus_{k\not=i,j}\hat\xi^{1}(x_{k})\oplus\hat\xi^1(y)\big)
\cap \big(\hat\xi^{1}(x_{i})\oplus \hat\xi^{1}(x_{j})\big)
=\lim_{m\rightarrow\infty} (F^m_{ij}(\dot z)).
$$
By Proposition \ref{hyperincrease}, all maps $F^m_{ij}$
are monotone. It follows that for all $t$ in $[\dot y, \dot z]$,
we have
$$
\lim_{m\rightarrow\infty} (F^m_{ij}(\dot y))=
\lim_{m\rightarrow\infty} (F^m_{ij}(t)).
$$
But for $t$ in a neighbourhood of $y$, the following sums
are direct
\begin{eqnarray*}
\forall i, \ \ \sum_{i\not=j}\hat\xi^1(x_j)+\xi^1(t)&=&E, 
\end{eqnarray*}
This implies that
$$
\lim_{m\rightarrow\infty} (F^m_{ij}(t)) =
\big(\bigoplus_{k\not=i,j}\hat\xi^{1}(x_{k})\oplus\hat\xi^1(\pi(t))\big)
\cap \big(\hat\xi^{1}(x_{i})\oplus \hat\xi^{1}(x_{j})\big)
$$
Hence, for all $t$ in a non empty open set, we have
$$
\hat\xi^1(t)=\hat\xi^1(y).
$$
But this is impossible, since $\hat\xi^{1}$ is
parametrised by arc length and, hence, not locally constant.

{\em $\Lambda$-Hyperconvexity: proof of Proposition \ref{bilhyp}.} In this proof, we shall only use the following property of the elements of $\Lambda$: If $U$ is a neighbourhood of an element of $\Lambda$, $\xi(U)$ is not included in a union of two proper subspaces (cf. Lemma \ref{localworm}$(iv)$). This property is then  true for any element in $\overline\Lambda$, the closure of $\Lambda$.

Let $p$ be the smallest integer, less than $n$, if it exists,  such that there exist $p$  points 
$(x_{1},\ldots,x_{p-2},y,z)$ cyclically oriented with $x_i,y,z,\in\overline\Lambda$ such that the following sum is
not direct 
$$
H=\sum_{i=1}^{i=p-2}\hat\xi^{1}(x_{i})+\hat\xi^{1}(y)+\hat\xi^{1}(z).
$$
Thanks to Proposition \label{bil}, $p\geq 3$.
Notice that, by minimality of $p$, the
following sums are direct and equal
\begin{eqnarray}
\sum_{i=1}^{i=p-2}\hat\xi^{1}(x_{i})+\hat\xi^{1}(y)
=\sum_{i=1}^{i=p-2}\hat\xi^{1}(x_{i})+\hat\xi^{1}(z)=H.
\label{hyperyz}
\end{eqnarray}
By our initial remark, we can now  choose $(x_{p-1},\ldots,x_{n})$  in an arbitrarily small neighbourhood $J$ of $x_1$ in $\overline\Lambda$,  such that the
following sums are direct
\begin{eqnarray}
\hat\xi^1(x_1)+\ldots+\hat\xi^1(x_n)&=&E,\cr
\forall i\geq p-1,\ \ \  \sum_{j\not=i}\hat\xi^1(x_j)+\hat\xi^1(y)&=&
E\label{goodsumdir2}\cr
\forall  i\geq p-1,\ \ \  \sum_{j\not=i}\hat\xi^1(x_j)+\hat\xi^1(z)&=&
E. 
\end{eqnarray}
Let $I$ be an interval containing $y$ and $z$ and none of the $x_i$.
Like in the previous proof, for any of the points $t\in\{x_1,\ldots,x_n,y,z\}$, we denote by $\dot t$ a
point such that $\pi (\dot t)=t$.  For any distinct integers
$i,j$, let's write
$$
W_{ij}=(\dot x_i,\dot x_j), \ \ Y_{ij}=(\ldots,\dot
x_l,\ldots)_{l\not\in\{i,j\}}.
$$
We can as in Section \ref{hyperincrease0} consider the
maps $F^m_{ij}$ defined for $i,j\geq p-1$, by
\begin{eqnarray*}
\mapping{ I}
{\mathbb P \big(\tilde\xi_{m}^{(2)}(W_{ij})\big)\setminus
\{\tilde\xi^1_m(\dot x_i)\}}
{t}
{\big(\tilde\xi_{m}^{(n-2)}(Y_{ij})\oplus\tilde\xi^1_{m}(t)
\big)\cap\tilde\xi_{m}^{(2)}(W_{ij})}.
\end{eqnarray*}
By Assertions (\ref{goodsumdir2}) and (\ref{hyperyz}), we
obtain that for all $i,j\geq p-1$,
$$
\lim_{m\rightarrow\infty} (F^m_{ij}(\dot y))
=\big(\bigoplus_{k\not=i,j}\hat\xi^{1}(x_{k})\oplus\hat\xi^1(y)\big)
\cap \big(\hat\xi^{1}(x_{i})\oplus \hat\xi^{1}(x_{j})\big)
=\lim_{m\rightarrow\infty} (F^m_{ij}(\dot z)).
$$
By Proposition \ref{hyperincrease}, all maps $F^m_{ij}$
are monotone. It follows that for all $t$ in $[\dot y, \dot z]$,
we have
$$
\lim_{m\rightarrow\infty} (F^m_{ij}(\dot y))=
\lim_{m\rightarrow\infty} (F^m_{ij}(t)).
$$
But for $t$ in a neighbourhood of $y$, the following sums
are direct
\begin{eqnarray*}
\forall i, \ \ \sum_{i\not=j}\hat\xi^1(x_j)+\xi^1(t)&=&E, 
\end{eqnarray*}
This implies that
$$
\lim_{m\rightarrow\infty} (F^m_{ij}(t)) =
\big(\bigoplus_{k\not=i,j}\hat\xi^{1}(x_{k})\oplus\hat\xi^1(\pi(t))\big)
\cap \big(\hat\xi^{1}(x_{i})\oplus \hat\xi^{1}(x_{j})\big)
$$
Hence, for all $t$ in a right neighbourhood $U^+$ of $y$, we have
$$
\bigoplus_{i=1}^{i=p-2}\hat\xi^{1}(x_{i})\oplus\hat\xi^{1}(t)
=\bigoplus_{i=1}^{i=p-2}\hat\xi^{1}(x_{i})\oplus\hat\xi^{1}(y)
$$
Hence, there exists a proper subspace $H^+$ of $E$, such that 
$$
\forall t\in U^+,\  \  \  \hat\xi^1(t)\subset H^+.
$$
A symmetric reasoning (after some cyclic permutation of $(x_1,\ldots,y,z)$),
 shows there exist a left neighbourhood $U^-$ of $y$, a proper subspace $H^-$ of $E$ such that
$$
\forall t\in U^-,\ \ \  \hat\xi^1(t)\subset H^-.
$$
As a conclusion we obtain a neighbourhood $U$ of $y$ and two proper subspaces $H^+$ and $H^-$ of $E$ such that
$$
\forall t\in U, \ \ \ \hat\xi^1(t)\subset H^-\cup H^+.
$$
This contradicts our initial remark.
\qed

\subsubsection{Conjugating to the action on the boundary at infinity}
We place ourselves in the situation described by Proposition
\ref{coi}. That is we consider
\begin{enumerate}
\item a sequence $\{\rho_{m}\}_{m\in\mathbb M}$ of
  representations of $\grf$ in $PSL(E)$, converging to $\rho$, such that for all
  non trivial $\gamma$ in $\grf$, $\rho_{m}(\gamma)$ is
  purely loxodromic \label{conjug1}
\item a sequence of maps $\{\xi^{1}_{m}\}_{m\in\mathbb N}$
  from $\partial_{\infty}\grf$ to $\mathbb P(E)$, such that
  each $\xi^{1}_{m}$ is $\rho_{m}$-equivariant,
  \label{conjug2}
\item a sequence of homeomorphisms $\phi_{m}$ of $S^{1}$ to
  $\partial_{\infty}\grf$ such that
  $\{\xi^{1}_{m}\circ\phi_{m}\}_{m\in\mathbb N}$ converges
  to $\hat\xi^{1}\circ\pi$ such that $\hat\xi^{1}$ is an
  embedding and $\pi$ is a monotone map from $S^{1}$ to
  $\partial_{\infty}\grf$. \label{conjug3}
\end{enumerate}

We say in this situation that $(\hat\xi^{1},\rho)$ is a {\em
  good limit}. Notice in this case, there exists an action
of $\grf$ on $S^{1}$ so that $\hat\xi^{1}$ is
$\rho$-equivariant.

The next lemma finishes the proof of Proposition \ref{lbi}.
\begin{lemma}\label{northsouth}
  Let $\grf$ be a surface group. Let $\rho$ be a
  representation of $\grf$ in $SL(E)$. Assume
\begin{itemize}
\item the restriction of $\rho$ to every finite index
  subgroup is irreducible,
\item there exists a map $\hat\xi^{1}$ from $S^{1}$ to
  $\mathbb P(E)$ such that $(\hat\xi^1,\rho)$ is a good
  limit (in particular every element of $\rho(\grf)$ has
  real eigenvalues).
\end{itemize}
Then the induced action of $\grf$ on $S^1$ for which
$\hat\xi^{1}$ is $\rho$-equivariant is topologically semi-conjugate to the action on $\partial_\infty \grf$ in the following sense : there exists an orientation preserving, left continuous,  $\grf$-equivariant map $\varphi_0$ from $\partial_\infty\grf$ to $S^1$.

Finally, $\hat\xi^1\circ\varphi_0$ is $*$-hyperconvex.
 \end{lemma}
 
 \proof We first notice that $\Gamma_{\mathbb R}$ is not empty. Indeed, the connected component of the  Zariski closure of $\rho(\grf)$ is irreducible. It therefore contains a non trivial diagonalisable element. Hence so does $\rho(\grf)$, but since this element has only real eigenvalues, this shows $\Gamma_{\mathbb R}$ is not empty. Let now $\Gamma^0$ be the subgroup of finite index of orientation preserving  elements of $\grf$ acting on $S^1$. Let (cf. Paragraph \ref{lamb}). Again $\Gamma^0_{\mathbb R}$, is not empty and invariant by conjugation.
 $$
 \Lambda=\Lambda_{\xi,\rho,\Gamma^0}\subset S^1.
 $$
 The set $\Lambda$ is not empty by Lemma \ref{localworm}.$(iii)$ and invariant under the action of $\grf$.
 
Let  $\gamma$ be an element of  $\Gamma^0_{\mathbb R}$,
\begin{itemize}
\item  let $\gamma^+$ be such that
$\xi(\gamma^+)$ is the attractive fixed point of $\rho(\gamma)$; the point $\gamma^+$  is well defined by Lemma \ref{localworm}.$(iii)$; let 
$$
\Lambda=\{\gamma^+, \gamma \in \Gamma^0_{\mathbb R}\},
$$
\item  Let
$\gamma_0^+$ be the attractive fixed
point of $\rho_0(\gamma)$, and
$$
\Lambda_0=\{\gamma_0^+, \gamma \in \Gamma^0_{\mathbb R}\},
$$
\end{itemize}

We know that
$$
\gamma_0^+=\lambda_0^+\\ \implies\\ \exists p,q\not=0,\\
\lambda^p=\gamma^q\\ \implies\\ \gamma^+=\lambda^+.
$$

Therefore we have a well defined map $\varphi_0$, maybe not
injective, defined from $\Lambda_0$ to
$\Lambda$. We notice that the first set is dense by the
minimality of the action of $\rho_0(\Gamma)$.

We now prove that $\varphi_0$ preserves the cyclic ordering. We
are going to use our full hypothesis concerning the
construction of $\hat\xi^1$, in particular that $\hat\xi^1$
is a good limit.  We begin by some observation.  By
construction, $\xi^{1}_{m}(\gamma_{0}^{+})$ is an attractive
fixed point of $\rho_{m}(\gamma)$ ({\em cf} Hypothesis
(\ref{conjug1}) and (\ref{conjug2})), hence
\begin{eqnarray}
\lim_{m\rightarrow\infty}(\xi^{1}_{m}(\gamma_{0}^{+}))=\hat\xi^1
\label{conjug4}
(\gamma^{+})
\end{eqnarray}
We can now extract a subsequence so that
$\{\tilde\gamma_{m}^{+}\}_{m\in\mathbb N}
=\{\phi_{m}^{-1}(\gamma_{0}^{+})\}_{m\in\mathbb N}$
converges to a point $\tilde\gamma^{+}$. By Hypothesis
(\ref{conjug3}),
\begin{eqnarray}
\lim_{m\rightarrow\infty}(\xi^{1}_{m}(\phi_{m}(\tilde\gamma_{m}^{+})))
=\hat\xi^1\circ\pi
(\tilde\gamma^{+})\label{conjug5}
\end{eqnarray}
Combining Assertions (\ref{conjug5}) and (\ref{conjug4}), and using the injectivity of $\hat\xi^1$ we
get that $\pi(\tilde\gamma^{+})=\gamma^{+}$. It follows that
$\varphi_{0}$ preserves the orientation. Indeed, let $\gamma$,
$\lambda$ and $\delta$ be three elements of $\Gamma$ such
that $(\gamma_{0}^{+}, \lambda_{0}^{+}, \delta_{0}^{+})$ are
cyclically oriented. Then so are $(\phi_{n}
(\gamma_{0}^{+}), \phi_{n} (\lambda_{0}^{+}), \phi_{n}
(\delta_{0}^{+}))$, hence $(\tilde\gamma_{0}^{+},
\tilde\lambda_{0}^{+}, \tilde\delta_{0}^{+})$ and finally
$(\pi(\tilde\gamma_{0}^{+}), \pi(\tilde\lambda_{0}^{+}),
\pi(\tilde\delta_{0}^{+}))=(\gamma^{+}, \lambda^{+},
\delta^{+})$.

Now that we know that $\phi_{0}$ preserves the orientation,
we can extend it by left continuity  to a $\Gamma$ equivariant orientation preserving map  from
$\partial_{\infty}\grf$ to $S^{1}$ since
$\Lambda_0$ is dense by the
minimality.

Let's finally prove $\varphi_{0}$ is injective. For that let
$$
U=\{x\in S^1/ \varphi_0 \hbox{ is constant on a
  neighbourhood of } x\}.
$$
Let's show that $U=\emptyset$. We can remark that $U$ is
an open $\rho_{0}(\grf)$-invariant strict subset of
$\partial_{\infty}\grf$. By minimality of the action of
$\grf$ on $\partial_\infty \grf$ we conclude that
$U=\emptyset$.  This imply that $\varphi_0$ is strictly
monotone, hence injective.

It remains to prove that $\hat\xi^1\circ\varphi_0$ is $*$-hyperconvex.
First, we notice that $\varphi_0(\partial_\infty\grf)=\overline\Lambda$.
By Proposition   \ref{bil}, we have  that for $n$ distinct points $(x_1,\ldots,x_n)$ of $\overline\Lambda$, the following sum is direct
$$
\sum_i \hat\xi^1(x_i).
$$
This implies the first condition on $*$-hyperconvexity.
The last condition on $*$-hyperconvexity is a closed condition and follows from the fact that $\xi^1\pi$ is a limit of $\{\xi_m^1\circ\phi_m\}$ which is a sequence of  hyperconvex maps.
\qed

\section{Appendix: some lemmas}\label{9}

We prove several lemmas that are used several times in the
article. This appendix is independent of what precedes. Next
lemma is obvious for Higgs field experts and certainly well
known.
\begin{lemma}\label{corohitchin}
  If $\rho$ belongs to Hitchin's component, then the
  connected component of the Zariski closure of $\rho(\grf)$
  is irreducible, or equivalently $\rho$ restricted to every
  finite index subgroup is irreducible.
\end{lemma}
\proof We have to recall a little bit of Hitchin's
construction. We consider a surface $S$, its canonical
bundle $K$ and the holomorphic vector bundle
$$
E=K^{-n \otimes}\oplus K^{(2-n)\otimes
}\oplus\ldots\ldots K^{n\otimes}.
$$
We consider the Higgs field $\phi$ which is a section of
$End(E)$ associated to a companion matrix.

Every representation in Hitchin's component comes from such
a Higgs field. From Lemma 1.2 in \cite{S}, the parallel
sections of the endomorphism bundle are exactly those
holomorphic sections commuting with the Higgs field. Let $A$
be such a section. Since $A$ is holomorphic the first row of
its matrix in the decomposition of $E$ vanishes. Now it is
easy to check that a matrix whose first row vanishes and
that commutes with a companion matrix is zero.

We have just proved that the endomorphism bundle has no
parallel sections. Hence the representation is irreducible.

Since the restriction to a finite index subgroup comes from
a representation in Hitchin's component on the corresponding
finite cover, we get the second part of the statement.  \qed

Next result is used several times.
\begin{lemma}\label{conseqirredu}
  Let $\Gamma$ be a surface group.  Let $\rho$ be a
  representation of $\grf$ to $SL(n,\mathbb R)$. Assume
  there exists a continuous $\rho$-equivariant map $\xi^{k+1}$ from
  $\partial_{\infty}\grf$ to the Grassmannian of $k+1$-planes
  in $\mathbb R^n$. Assume there exists some $(n-k)$-plane
  A, a non empty open set $U$ in $\partial_{\infty}\grf$
  such that
\begin{equation}
\forall y\in U,\ \dim(\xi^{k+1}(y)\cap A)\geq 1\label{qfirred1}.
\end{equation}
Then, the restriction of $\rho$ to a finite index subgroup
is not irreducible, or equivalently the connected component
of the identity of the Zariski closure of $\rho(\grf)$ is
not irreducible.
\end{lemma} 

\proof {\it First step.} We show first that there exists
some $(n-k)$-plane $B$, such that
\begin{equation}
\forall y\in \partial_{\infty}\grf,\ \dim(\xi^{k}(y)\cap B)\geq 1\label{qfh2}.
\end{equation}
First, we can find a smaller open subset $O$ of $U$ and some
element $\gamma\in\grf$ such that
\begin{eqnarray*}
\gamma^{i}(O)&\subset&\gamma^{i+1}(O) \\
O_\infty&=&\bigcup_{i\in\mathbb N}\gamma^i(O)\hbox{ is dense in }\partial_{\infty}\grf.
\end{eqnarray*}
Next notice that for $B^i=\rho(\gamma)^{-i}(A)$ we have
\begin{equation*}
\forall y\in\gamma^i(O) ,\ \dim(\xi^{k}(y)\cap B^i)\geq 1.
\end{equation*}
Now we extract from $\{B^i\}_{i\in\mathbb N}$ a convergent
subsequence to a $(n-k)$-plane $B$. It follows that
\begin{equation*}
\forall y\in O_\infty,\ \dim(\xi^{k}(y)\cap B)\geq 1.
\end{equation*}
Hence, by density of $O_\infty$, we obtain Assertion
(\ref{qfh2}).  The proof finally follows from the next
proposition applied to $G$ the Zariski closure of
$\rho(\grf)$. \qed

\begin{proposition}\label{algebirredu}
  Let $G$ be an algebraic subgroup of $SL(n,\mathbb R)$. If
  there exist a $k$-plane $C$, a $(n-k)$-plane $B$ such that
\begin{equation}
\forall g\in G,\ \dim(g(C)\cap B)\geq 1\label{qfh3}.
\end{equation}
Then the connected component of the identity of $G$ is not
irreducible.
\end{proposition}
\proof Let $G$, $C$, $B$ be as in the Proposition.  Notice
that $B+C$ is a proper subspace, and $B\cap C$ is not
reduced to $\{0\}$.  Let $g_0$ be an element of $G$ such
that
$$
\forall g\in G,\ p:=\dim(g_0(C)\cap B)\leq\dim(g(C)\cap
B).
$$
Let now $D=g_0(C)$, $F$ be a codimension $p-1$ vector
subspace of $B$ such that,
$$
\dim(D\cap F)=1.
$$
Notice now that
\begin{equation}
\forall g\in G,\ \dim(g(D)\cap F)\geq 1\label{qfh4}.
\end{equation}
Let $\mathcal G$ be the Lie algebra of $G$.

We will show the following assertion that contradicts the
irreducibility of the connected component of the identity of
$G$,
\begin{equation}
\forall \alpha\in \mathcal G,\ \alpha(D\cap F) \subset D+F.\label{qfh5}
\end{equation}
Indeed, let $(e_0,\ldots,e_k)$ be a basis of $D$, let
$(u_1,\ldots,u_l)$ be a basis of $F$, such that $(u_l)$ is a
basis of $D\cap F$. Notice that
$$
e_0\wedge\ldots\wedge e_k\wedge u_1\wedge\ldots\wedge
u_{l-1}\not=0.
$$
Let $h$ be an element of $\mathcal G$. From Assertion
(\ref{qfh4}), we obtain
$$
e_0\wedge \ldots\wedge e_k\wedge
e^{t\alpha}(u_1\wedge\ldots\wedge u_{l})=0.
$$
We take now the first order term of the above series in
$t$. Since
$$
e_0\wedge \ldots\wedge e_k\wedge u_l=0,
$$
we obtain,
$$
e_0\wedge \ldots\wedge e_k\wedge u_1\wedge\ldots\wedge
u_{l-1}\wedge\alpha (u_{l})=0.
$$
This implies Assertion (\ref{qfh5}).\qed

Finally, the next lemma is of independent interest

\begin{lemma}\label{discrete}
  Let $\Gamma$ be a subgroup of $SL(n,\mathbb R)$ whose
  elements are all real split. Assume every finite index
  subgroup of $\Gamma$ is irreducible. Then $\Gamma$ is
  discrete.
\end{lemma}
\proof Let $G$ be the Zariski closure of $\Gamma$.  From the
irreducibility hypothesis it follows $G$ is semi-simple. If
$\Gamma$ is not discrete, since it is Zariski dense, it is
classical that its closure (for the usual topology) contains
one of the non trivial factor $H$ of $G$.  But the closure
of $\Gamma$ consists of elements whose elements have only
real eigenvalues.  This would imply the maximal compact
subgroup of $H$ is reduced to the identity and this never
happens for a simple Lie group.  \qed

\auteur

\begin{thebibliography}{99}
\bibitem{IND}{I. Biswas, P. ArŽs-Gastesi, S. Govindarajan
  } {\em Parabolic Higgs bundles and Teichmüller spaces
    for punctured surfaces.}  {Trans. Amer. Math. Soc. {\bf
      349} (1997), no. 4, 1551--1560.}
\bibitem{BIW}{M, Burger, A. Iozzi, A. Wienhard }{\em Sur les représentations d'un groupe de surface compacte avec invariant de Toledo maximal.}{ C.R. Acad. Sci. Paris, Ser. I 336 (2003), 387-390}
\bibitem{CG}{S. Choi, W. Goldman }{\em Convex real
    projective structures on closed surfaces are closed.}{
    Proc. Amer. Math. Soc.  {\bf 118} (1993), no. 2,
    657--661.}
\bibitem{KC} {K. Corlette }{\em Flat G-bundles with
    canonical metrics.}{ J. Differential Geom.  {\bf 28}
    (1988), no. 3, 361--382.}
\bibitem{D}{S.K. Donaldson }{\em Twisted harmonic maps and
    the self-duality equations.}{ Proc. London Math. Soc.
    (3) {\bf 55} (1987), 127-131}
\bibitem{FLP}{A. Fathi, F. Laudenbach, V. Poenaru} {\em
    Travaux de Thurston sur les surfaces}{ Astérisque,
    No.66-67. Paris: Société Mathématique de
    France. 284 p. (1979).}
    \bibitem{GB}{V. Fock, A.B Goncharov} {\em  Moduli spaces of local systems and higher Teichmuller theory
,}{math.AG/031114  }
\bibitem{G1}{W. Goldman }{Topological components of spaces
    of representations.}{Invent. Math. {\bf 93} (1988), no.
    3, 557--607.}
    \bibitem{OG}{O. Guichard} {\em Sur les répresentations des groupes de surface} {preprint}
       \bibitem{WG}{W. Goldman, }{\em Geometric Structures and varieties of representation}{in ``The Geometry of Group Representations"}{Summer Conference 1987, Bouldern Colorado, W. Goldman and A. Magid (eds)} 
\bibitem{UH}{U. Hamenstädt} {\em Cocycles, Hausdorff measures and cross ratios }{ Ergodic Theory and Dynamical Systems (1997), 17 pp 1061-1081}
\bibitem{HK}{B. Hasselblatt, A.Katok }{\em Introduction to
    the modern theory of dynamical systems. With a
    supplementary chapter by Katok and Leonardo Mendoza. } {
    Encyclopedia of Mathematics and its Applications, 54.
    Cambridge University Press, Cambridge, 1995.}
    \bibitem{H}{N. Hitchin }{\em Lie Groups and Teichmüller
    spaces.}{ Topology {\bf 31} (1992), no. 3, 449--473.}
\bibitem{MK}{M. Kontsevich }{\em The Virasoro algebra and
    Teichmüller spaces.  (Russian)}{ Funktsional. Anal. i
    Prilozhen. 21 (1987), no. 2, 78--79}

\bibitem{FL1}{F. Labourie }{\em Existence d'applications
    harmoniques tordues à valeurs dans les variétés
    à courbure négative.}{ Proc. Amer. Math. Soc.  111
    (1991), no. 3, 877--882.}
\bibitem{FL2}{F. Labourie }{\em $\mathbb{RP}^2$-structures
    et differentielles cubiques holomorphes.}{ Proceedings
    of the GARC Conference in Differential Geometry, Seoul
    National University, 1997. }
\bibitem{FL3}{F. Labourie }{\em Crossratios, Surface Groups
    and $SL(n,\mathbb R)$.}{ preprint}
 \bibitem{FL4}{F. Labourie }{\em Crossratios, Anosov Representations and Quasi-Isometries}{in preparation}
   
\bibitem{led}{F. Ledrappier }{\em Structure au bord des
    variétés à courbure négative. }{
    Sé\-mi\-nai\-re de Théorie Spectrale et
    Géométrie, No. 13, Année 1994--1995, 97--122,
    Sémin. Théor. Spectr. Géom., 13, Univ.
    Grenoble I, Saint-Martin-d'Hères, 1995. }
    
\bibitem{JL}{J. Loftin }{\em Affine spheres and convex
    $\mathbb{RP}\sp n$-manifolds.}{ Amer. J. Math.  {\bf
      123} (2001), no. 2, 255--274.}

    \bibitem{LO}{W. L. Lok, }{\em Deformation of locally homogeneous spaces and Kleinian groups }{ Doctoral Thesis, Columbia University 1984}
\bibitem{JPO}{J.-P. Otal }{\em Le spectre marquéŽ des
    longueurs des surfaces à courbure négative.}{Ann.
    of Math. (2) {\bf 131} (1990), no. 1, 151--162.}
\bibitem{AP}{A. Parreau }{\em Dégénérescences de
    sous-groupes discrets de groupes de Lie semisimples et
    actions de groupes sur des immeubles affines }{ Thèse
    Université Paris-Sud (2000)}
\bibitem{P}{R. Penner }{\em The decorated Teichmüller
    space of punctured surfaces.}{ Comm. Math. Phys. 113
    (1987), no. 2, 299--339.}
\bibitem{CR}{C. Robinson }{\em Stability, Symbolic Dynamics,
    and Chaos.  }{ Studies in Advanced Mathematics.  CRC
    Press, Boca Raton, FL, 1995. }
\bibitem{GS}{G. Segal }{ \em The geometry of the KdV
    equation.}{ Topological methods in quantum field theory
    (Trieste, 1990).  Internat. J. Modern Phys. A 6 (1991),
    no. 16, 2859--2869.}
\bibitem{S}{C. Simpson}{ \em Higgs bundles and local
    systems. }{ Inst. Hautes Études Sci. Publ. Math.  No.
    75, (1992), 5--95.}




\end{thebibliography}
\end{document}